\documentclass[a4paper,oneside,10.5pt]{article}%
\usepackage{amsmath}
\usepackage{amsfonts}
\usepackage{indentfirst}
\usepackage{amssymb}
\usepackage{graphicx}%
\usepackage{appendix}
\usepackage[colorlinks,citecolor=blue,urlcolor=blue]{hyperref}

\providecommand{\U}[1]{\protect \rule{.1in}{.1in}}

\newtheorem{theorem}{Theorem}[section]

\newtheorem{assumption}[theorem]{Assumption}
\newtheorem{example}[theorem]{Example}

\newtheorem{lemma}[theorem]{Lemma}

\newtheorem{remark}[theorem]{Remark}

\numberwithin{equation}{section}

\begin{document}

\title{A varying terminal time structure for stochastic optimal control under constrained condition}
\author{Shuzhen Yang\thanks{Shandong University-Zhong Tai Securities Institute for Financial Studies, Shandong University, PR China, (yangsz@sdu.edu.cn).}
\thanks{This work was supported by the National Natural Science Foundation of China (Grant No.11701330, 11871050) and Young Scholars Program of Shandong University.}}
\date{}
\maketitle

\textbf{Abstract}:   In this study, we propose a varying terminal time structure for the optimal control problem under state constraints, in which the terminal time follows the varying of the control via the constrained condition. Focusing on this new optimal control problem, we investigate a novel stochastic maximum principle, which differs from the traditional optimal control problem under state constraints. The optimal pair of the optimal control model can be verified via this new stochastic maximum principle.

\textbf{Keywords}: varying terminal time; stochastic differential equation; stochastic maximum principle

{\textbf{MSC2010}: 93E03; 93E20; 60G99

\addcontentsline{toc}{section}{\hspace*{1.8em}Abstract}

\section{Introduction}
In  the traditional stochastic optimal control problem, we usually consider the following model. For a given positive constant $T$, we denote the running cost by $f(X(t),u(t))$ at time $t\in[0,T]$, and the terminal cost at time $T$ is given  by $\Psi(X(T))$. The cost functional is given as follows:
\begin{equation}
\label{incos-1}
J(u(\cdot))=
\mathbb{E}\bigg{[}\displaystyle\int_0^Tf(X^u(t),u(t))\mathrm{d}t+\Psi(X^u(T))\bigg{]},
\end{equation}
where the state process $X^u(\cdot)$ is driven by the following controlled stochastic differential equation:
\begin{equation}
\label{insde-1}
X^u(s)=x_0+\int_{0}^{s}b(X^u(t),u(t))\mathrm{d}t+\int_{0}^{s}\sigma(X^u(t),u(t))\mathrm{d}W(t),\ s\in [0,T].
\end{equation}

The stochastic maximum principle and dynamic programming principle constitute two powerful tools for studying the traditional stochastic optimal control problem. We refer the reader to Bensoussan \cite{B81} and Bismut \cite{B78} concerning the local maximum principle with a convex control set, and to Peng \cite{P90}  for the global  maximum principle with a general control domain. Furthermore, We refer the reader to Hu \cite{17H} for the stochastic global maximum principle for recursive utilities systems, among others \cite{HJX18,W13,Y10}. In addition, Hu and Ji \cite{16HJ}  studied the  stochastic maximum principle for the stochastic recursive optimal control problem under volatility ambiguity.  L\"u and Zhang \cite{LZ14} studied the general stochastic maximum principle for a backward stochastic evolution equation in infinite dimensions. Qiu and Tang \cite{QT12} studied the maximum principle for quasi-linear backward stochastic partial differential equations. Yang \cite{Y18} studied  stochastic differential systems with a multi-time states cost functional. Further details regarding  theories of optimal control problems can be found in  \cite{FM06,Y99}.

For the optimal control problem under state constraints, Frankowska \cite{F10} reviewed the basic theory of deterministic optimal controls, and the necessary optimality conditions under state constrained.  Rutquist\cite{R17} reviewed a few methods to solve stochastic optimal control problems under state constraints. Bouchard et al. \cite{BE10} considered the stochastic optimal control problems in which the controlled process satisfies an almost sure constraint at final time. Furthermore,  L. Bourdin and E. Tr\'{e}lat studied  more general versions of the Pontryagin maximum principle on time scales in \cite{BT13,BT17}.

 In the traditional optimal control problem, cost functional (\ref{incos-1}) is minimized within a given length of time $T$ by the control state process (\ref{insde-1}). In general, it is useful not to specify the terminal time before beginning to control the state process (\ref{insde-1}), which allows that the terminal time can depend on the values of $X^u(\cdot)$ or the control $u(\cdot)$. For example, in the investment portfolio problem, the state process  (\ref{insde-1}) can be used to describe the asset, while cost functional (\ref{incos-1}) can be used to represent the risk of the asset $X^u(\cdot)$. First, we specify a positive constant $T$ that denotes the period of investment. We can stop the investment plan before the period $T$. The criterion for stopping the investment can be described as follows:
\begin{equation}
\label{intime-1}
\tau^u=\inf\bigg{\{}t:\mathbb{E}[X^u(t)]\geq \alpha,\ t\in [0,T] \bigg{\}}\bigwedge T,
\end{equation}
{ where $a_1\bigwedge a_2=\min(a_1,a_2),\ a_1,a_2\in \mathbb{R}$} and $\alpha$ describes the target  of the mean value of the asset $X^u(\cdot)$ in the period $T$. Thus we only need to minimize the risk within $[0,\tau^u]$, and the cost functional is given as follows:
\begin{equation}
\label{incos-2}
J(u(\cdot))=\mathbb{E}\bigg{[}\displaystyle\int_0^{\tau^u}f(X^u(t),u(t))\mathrm{d}t+\Psi(X^u(\tau^u))\bigg{]}.
\end{equation}
As shown in (\ref{intime-1}), unlike that in the traditional optimal control problem, the terminal time $\tau^u$ varies according to the control $u(\cdot)$. Specially, in the investment  portfolio problem, we need to balance the mean value and risk, or equivalently, to balance the terminal time $\tau^u$ and  cost functional $J(u(\cdot))$. We considered this varying terminal time mean-variance problem in Yang \cite{Y19}, in which an optimal strategy  and the related varying terminal time are found.

{ Considering the optimal control problem under state constraints $\mathbb{E}[\Phi(X^u(T))]\geq \alpha$, we introduce a varying terminal time optimal control structure:
\begin{equation}
\label{intime-2}
\tau^u=\inf\bigg{\{}t:\mathbb{E}[\Phi(X^u(t))]\geq \alpha,\ t\in [0,T] \bigg{\}}\bigwedge T.
\end{equation}
Notice that in (\ref{intime-2}), we consider the minimum time $\tau^u$ such that
$\mathbb{E}[\Phi(X^u(\tau^u))]\geq \alpha,\ \tau^u\in [0,T]$ with a control $u(\cdot)$. If $\bigg{\{}t:\mathbb{E}[\Phi(X^u(t))]\geq \alpha,\ t\in [0,T] \bigg{\}}\neq \varnothing$, the constrained condition in $\tau^{u}$ can be viewed as the state constraints
$$
\mathbb{E}[\Phi(X^u(t))]< \alpha, \ t\in [0,\tau^{u}),\  \mathbb{E}[\Phi(X^u(\tau^{u}))]\geq \alpha.
$$}
Focusing on this optimal control problem, we consider three  cases of $\tau^{\bar{u}}$ for $(\bar{u}(\cdot),\bar{X}(\cdot))$, where $(\bar{u}(\cdot),\bar{X}(\cdot))$ is a given optimal pair for the cost functional (\ref{incos-2}). Based on case $(i),\ \tau^{\bar{u}}<T$, case $(ii),\ \inf\bigg{\{}t:\mathbb{E}[\Phi(\bar{X}(t))]\geq \alpha,\ t\in [0,T] \bigg{\}}=T,$ and case $(iii),\ \bigg{\{}t:\mathbb{E}[\Phi(\bar{X}(t))]\geq \alpha,\ t\in [0,T] \bigg{\}}=\varnothing$. We calculate the variation of the varying terminal time $\tau^{\bar{u}}$ and achieve a novel stochastic maximum principle. Furthermore, we compare our optimal control problem with the traditional optimal control problem under state constraints.

The remainder of this paper is organized as follows: In Section 2, we formulate a varying terminal time stochastic optimal control problem. Then, we establish the stochastic maximum principle in Section 3, and compare our optimal control problem with the traditional optimal control problem under state constraints. Finally, we conclude the results of the paper and describe some possibilities future work in Section 4.

\section{The new optimal control problem}

Let $W$ be a $d$-dimensional standard Brownian motion defined on a complete
filtered probability space $(\Omega,\mathcal{F},P;\{ \mathcal{F}(t)\}_{t\geq
0})$, where $\{ \mathcal{F}(t)\}_{t\geq0}$ is the $P$-augmentation of the
natural filtration generated by the Brownian motion $W$. Let $T>0$ be given, considering the following controlled stochastic differential equation,
\begin{equation}
\text{d}{X}^u(t)=b(X^u{(t)},u(t))\text{d}t+\sigma (X^u{(t)},u(t))\text{d}W(t) ,\quad t\in(0,T],\label{ODE_1}%
\end{equation}
with the initial condition $X(0)=x_0$, { where
$u(\cdot)\in L^2_{\mathcal{F}}(0,T;U)$ is a control process taking value in a convex set $U$ of $\mathbb{R}^k$ with  a given positive integer $k$, and $L^2_{\mathcal{F}}(0,T;U)$ is the set of all $U$ valued, measurable processes $\phi(\cdot)$ adapted to $\{\mathcal{F}_t\}_{t\geq 0}$ such that
$$
\mathbb{E}\bigg{[}\int_0^{T}\left|\phi(t) \right|^2\mathrm{d}t \bigg{]}<+\infty.
$$}
In this study, we consider the following varying terminal time cost functional:
\begin{equation}
J(u(\cdot))=%
\mathbb{E}\bigg{[}{\displaystyle \int \limits_{0}^{\tau^{u}}}
f(X^u{(t)},u(t))\text{d}t+\Psi(X^u(\tau^u))\bigg{]},\label{cost-1}%
\end{equation}
where
\begin{equation}
\label{time-1}
\tau^u=\inf\bigg{\{}t:\mathbb{E}[\Phi(X^u(t))]\geq \alpha,\ t\in [0,T] \bigg{\}}\bigwedge T,
\end{equation}
$a_1\bigwedge a_2=\min(a_1,a_2),\ a_1,a_2\in \mathbb{R}$ and $\alpha\in (\Phi(x_0),+\infty)$. Note that if $\alpha<\Phi(x_0)$, then $\tau^u=0$, and the problem is trivial.  Furthermore,
\[%
\begin{array}
[c]{l}%
b:\mathbb{R}^m\times U\to \mathbb{R}^m,\\
\sigma:\mathbb{R}^{m}\times U\to \mathbb{R}^{m\times d},\\
f:\mathbb{R}^m\times U\to \mathbb{R},\\
\Psi,\Phi:\mathbb{R}^{m}\to \mathbb{R},\\
\end{array}
\]
we set $\sigma=(\sigma^1,\sigma^2,\cdots,\sigma^d)$, and $\sigma^j\in \mathbb{R}^m$ for $j=1,2,\cdots, d$. Here,
$\mathbb{R}^m=\mathbb{R}^{m\times 1}$. In addition,  "$\top$" denotes the transform of vector or matrix.

We assume that $b,\sigma,f$ are uniformly continuous and satisfy the following
linear growth and Lipschitz conditions.

\begin{assumption}
\label{ass-b}There exists a constant $c>0$ such that%
\[%
\begin{array}
[c]{c}%
\left| b(x_{1},u)-b(x_{2},u)\right| +\left| \sigma(x_{1},u)-\sigma(x_{2},u)\right|
 \leq c\left|x_1-x_2 \right|,\\
\end{array}
\]
$\forall(x_{1},u),(x_{2},u)\in{\mathbb{R}^m}\times U$.
\end{assumption}

\begin{assumption}
\label{assb-b2} There exists a constant $c>0$ such that
\[
\left|b(x,u)\right|+\left|\sigma(x,u)\right| \leq c(1+\mid x \mid),\quad \forall(x,u)\in{\mathbb{R}^m}\times U.
\]
\end{assumption}

\begin{assumption}
\label{ass-fai}Let $b,\sigma,f$ be differentiable at $(x,u)$, and their derivatives in $x$ be uniformly continuous in $(x,u)$. Let $\Psi$ be twice differentiable at $x$, with its  derivatives in $x$ be  uniformly continuous in $x$. Let $\Phi$ be three-times differentiable at $x$, and its derivatives in $x$ be uniformly continuous in $x$.
\end{assumption}

\begin{remark}
We assume strong smoothness conditions on $\Psi,\Phi$ in Assumption \ref{ass-fai}. This is because, we need to calculate the variation of the varying terminal time $\tau^{u}$. See Section \ref{smp0} and Appendix \ref{app-a} for further details.
\end{remark}

Let $\mathcal{U}[0,T]=L^2_{\mathcal{F}}(0,T;U).$ If
Assumptions \ref{ass-b} and \ref{assb-b2} hold, then there exists a unique
solution $X^u(\cdot)$ for equation (\ref{ODE_1}) (see \cite{LS78}). $\bar{u}(\cdot)\in \mathcal{U}[0,\tau^{\bar{u}}]$
satisfying
\begin{equation}
J(\bar{u}(\cdot))= \underset{u(\cdot)\in\mathcal{U}[0,\tau^u]}{\inf}J(u(\cdot)) \label{cost-2}%
\end{equation}
is called an optimal control. The corresponding state trajectory $(\bar{u}(\cdot),\bar{X}(\cdot))$ is called an optimal state trajectory or optimal pair and $\tau^{\bar{u}}$ is  called the optimal terminal time.

\section{Preliminary results and maximum principle}\label{smp0}
{ To investigate the well-known Pontryagin stochastic maximum principle for the varying terminal time optimal control problem, we show the preliminary results in Subsection \ref{pre} and stochastic maximum principle in Subsection \ref{smp}. The relationship between the main results of Subsection \ref{pre} and Subsection \ref{smp} is given in Figure 1. As shown in Figure 1: the formula of constrained condition function $\Phi(X^{u}(\cdot))$ is given in Lemma \ref{le-00}; the asymptotic behavior of optimal terminal time $\tau^{\bar{u}}$ and  optimal solution $\bar{X}(\cdot)$ are showed in Lemma \ref{le-0} and Lemma \ref{le-1}, respectively; the variational equation for optimal terminal time $\tau^{\bar{u}}$ and cost functional $J(\bar{u}(\cdot))$ are presented in Lemma \ref{le-2} and Lemma \ref{le-3}, respectively. Furthermore, we establish the stochastic maximum principle for the cost functional $J(\bar{u}(\cdot))$ with varying terminal time $\tau^{\bar{u}}$ in Theorem \ref{Max}.}
\begin{figure}[h]
\begin{center}
\includegraphics[width=5.0 in]{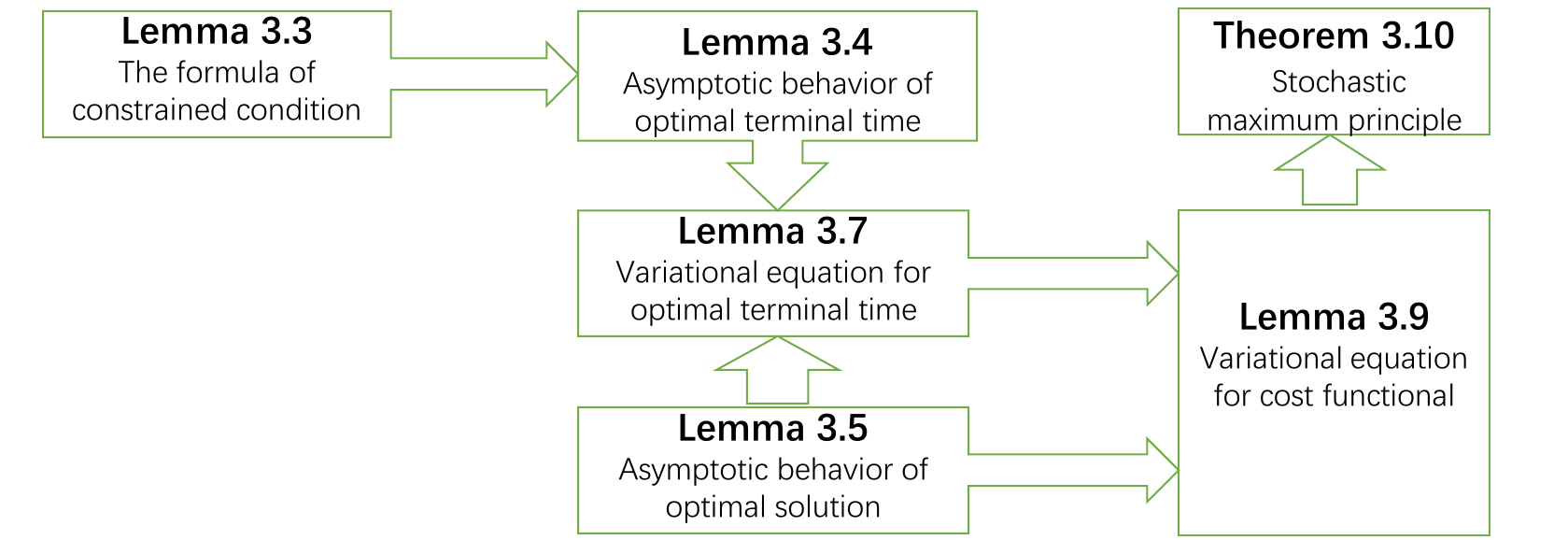}%
 \caption{The road map of the main results of Section \ref{smp0}}
\end{center}
\label{fig0}
\end{figure}

\subsection{Preliminary lemmas}\label{pre}
Notice that,
 $$
 \tau^{\bar{u}}=\inf\bigg{\{}t:\mathbb{E}[\Phi(\bar{X}(t))]\geq \alpha,\ t\in [0,T] \bigg{\}}\bigwedge T,
$$
and $\inf \varnothing=+\infty$.  Denoting
$$
\alpha_{\min}=\inf_{t\in [0,\tau^{\bar{u}}]}\mathbb{E}[\Phi(\bar{X}(t))],\quad \alpha_{\max}=\sup_{t\in [0,\tau^{\bar{u}}]}\mathbb{E}[\Phi(\bar{X}(t))],
$$
if $\alpha\in [\alpha_{\min},\alpha_{\max}]$, we have
$$
 \tau^{\bar{u}}=\inf\bigg{\{}t:\mathbb{E}[\Phi(\bar{X}(t))]\geq \alpha,\ t\in [0,T] \bigg{\}}.
$$

\begin{remark}
\label{re-1}
In this new stochastic optimal control problem, the terminal time $\tau^u$ varies with the value of $u(\cdot)$.
If we take $\alpha>\sup_{u\in \mathcal{U}[0,T]} \mathbb{E}[\Phi(X^{{u}}(t))]$, then $\tau^{u}=T$, and this new  optimal control problem  reduces to the traditional optimal control problem. If $\bigg{\{}t:\mathbb{E}[\Phi(X^u(t))]\geq \alpha,\ t\in [0,T] \bigg{\}}\neq \varnothing$, then combining the facts that $\alpha \in (\Phi(x_0),+\infty)$ and $\mathbb{E}[\Phi(X^u(\cdot))]$ is continuous in $t$, we have $ \mathbb{E}[\Phi(X^{{u}}(\tau^{{u}}))]=\alpha$.
\end{remark}

In the cost functional (\ref{cost-1}),  we consider a varying terminal time  cost functional, which is  different from that in the traditional optimal control problem. Note that $U$ is a convex set. Let $(\bar{u}(\cdot),\bar{X}(\cdot))$ be a given optimal pair, $0<\rho<1$,   and $v(\cdot)+\bar{u}(\cdot)\in \mathcal{U}[0,T]$ be any given control. We define
\[
u_{\rho}(t)=\bar{u}(t)+\rho v(t)=(1-\rho)\bar{u}(t)+\rho(v(t)+\bar{u}(t)),\ t\in [0,T].
\]
Clearly, $u_{\rho}(\cdot)\in \mathcal{U}[0,T]$, and
$X^{\rho}(\cdot)$ is the solution of equation (\ref{ODE_1}) under the
control $u_{\rho}(\cdot)$.
\begin{remark}
Note that we minimize the cost functional (\ref{cost-1}) on $[0,\tau^{u}]$, where $\tau^{{u}}$ depends on the expectation of $\Phi(X^u(\cdot))$, which implies that $\tau^{u_{\rho}}$ may be larger or smaller than $\tau^{\bar{u}}$. To give the definition of $u_{\rho}(\cdot)$ on $[0,\tau^{u_{\rho}}]$, we need to consider an optimal pair $(\bar{u}(\cdot),\bar{X}(\cdot))$ on $[0,T]$, where $\max(\tau^{\bar{u}},\tau^{u_{\rho}})\leq T$. In the following, we prove that $\left|\tau^{\bar{u}}-\tau^{u_{\rho}} \right|$ converges to $0$ as $\rho\to 0$ under certain continuity conditions.

\end{remark}

We first show the following result for the constrained condition function $\Phi(X^u(\cdot))$.
\begin{lemma}\label{le-00}
Let Assumptions \ref{ass-b}, \ref{assb-b2} and \ref{ass-fai} hold. We have
\begin{equation}
\label{mean-0}
\mathbb{E}[\Phi(X^u(s))]=\Phi(x_0)+\int_0^sh^u(t)\mathrm{d}t,\ s\in [0,T],
\end{equation}
where
$h^{{u}}(t)=\mathbb{E}\bigg{[}\Phi_x(X^{{u}}(t))^{\top}b(X^{{u}}(t),{u}(t))
+\displaystyle\frac{1}{2}\sum_{j=1}^d\sigma^j(X^{{u}}(t),{u}(t))^{\top}\Phi_{xx}(X^{{u}}(t))\sigma^j(X^{{u}}(t)
,{u}(t))\bigg{]}$, $t\in [0,s]$.

\end{lemma}
\textbf{Proof:} Applying  It\^{o} formula to $\Phi(X^u(t)),\ t\in [0,T]$,  we have
\begin{equation}%
\begin{array}
[c]{ll}%
& \text{d}\Phi(X^u(t))=\Phi_x(X^u(t))^{\top}\text{d}X^u(t)+\displaystyle\frac{1}{2}\big{[}\text{d}X^u(t)\big{]}^{\top}\Phi_{xx}(X^u(t))\text{d}X^u(t).\\
\end{array}
\end{equation}
From equation (\ref{ODE_1}), it follows that
\begin{equation*}%
\begin{array}
[c]{ll}%
 \text{d}\Phi(X^u(t))=&\bigg{[}\Phi_x(X^u(t))^{\top}b(X^u(t),u(t))+\displaystyle\frac{1}{2}\sum_{j=1}^d\sigma^j(X^u(t),u(t))^{\top}\Phi_{xx}(X^u(t))\sigma^j(X^u(t),u(t))\bigg{]}\text{d}t\\
&+\Phi_x(X^u(t))^{\top}\sigma(X^u(t),u(t))\text{d}W(t).\\
\end{array}
\end{equation*}
Thus, for $s\in [0,T]$,
\begin{equation*}%
\begin{array}
[c]{ll}%
 &\Phi(X^u(s))-\Phi(x_0)\\
=&\displaystyle \int_0^s\bigg{[}\Phi_x(X^u(t))^{\top}b(X^u(t),u(t))+\frac{1}{2}\sum_{j=1}^d\sigma^j(X^u(t),u(t))^{\top}\Phi_{xx}(X^u(t))\sigma^j(X^u(t),u(t))\bigg{]}\text{d}t\\
&+\displaystyle\int_0^s\Phi_x(X^u(t))^{\top}\sigma(X^u(t),u(t))\text{d}W(t).\\
\end{array}
\end{equation*}
Taking the expectation on both sides of the above equation, we have
\begin{equation*}%
\begin{array}
[c]{ll}%
 &\mathbb{E}[\Phi(X^u(s))]-\Phi(x_0)\\
=&\displaystyle \int_0^s\mathbb{E}\bigg{[}\Phi_x(X^u(t))^{\top}b(X^u(t),u(t))
+\frac{1}{2}\sum_{j=1}^d\sigma^j(X^u(t),u(t))^{\top}\Phi_{xx}(X^u(t))\sigma^j(X^u(t),u(t))\bigg{]}\text{d}t.\\
\end{array}
\end{equation*}
\noindent This completes the proof. $\ \ \ \ \ \ \ \ \Box$

{ Note that, the optimal varying terminal time is defined as
$$
 \tau^{\bar{u}}=\inf\bigg{\{}t:\mathbb{E}[\Phi(\bar{X}(t))]\geq \alpha,\ t\in [0,T] \bigg{\}}\bigwedge T.
$$
Denoting $A^{\bar{u}}=\bigg{\{}t:\mathbb{E}[\Phi(\bar{X}(t))]\geq \alpha,\ t\in [0,T] \bigg{\}}$, thus, $\tau^{\bar{u}}=\min(\inf A^{\bar{u}}, T)$. To calculate the variation of the varying terminal time $\tau^{\bar{u}}$, we introduce $u_{\rho}(\cdot)\in \mathcal{U}[0,T]$, where $ u_{\rho}(t)=\bar{u}(t)+\rho v(t),\ t\in [0,T]$.
To prove that $\tau^{{u}_{\rho}}$ converges to $\tau^{\bar{u}}$ when $\rho\to 0$,
we assume that the condition $H_0$ is right, where {$H_0$}: $h^{\bar{u}}(\tau^{\bar{u}})\neq 0$ and $h^{\bar{u}}(\cdot)$ is continuous at the point $\tau^{\bar{u}}$. Note that,
\begin{equation*}
\mathbb{E}[\Phi(X^{\bar{u}}(s))]=\Phi(x_0)+\int_0^sh^{\bar{u}}(t)\mathrm{d}t,\ s\in [0,T],
\end{equation*}
the condition $H_0$ can be used to guarantee that the constrained condition function $\mathbb{E}[\Phi(X^{\bar{u}}(\cdot))]$ is differentiable at $\tau^{\bar{u}}$ and the related derivative is continuous at $\tau^{\bar{u}}$ and does not equal to $0$. This is the key point to obtain the variation of the optimal varying terminal time $\tau^{\bar{u}}$.  Furthermore, we show the continuity of  $\tau^{u}$ at $\bar{u}(\cdot)$ in Lemma \ref{le-0}, and prove that $\tau^{u}$ is differentiable and continuous at $\bar{u}(\cdot)$ in Lemma \ref{le-2}.

Note that, $\tau^{\bar{u}}=\min(\inf A^{\bar{u}}, T)$, we have three cases about $\tau^{\bar{u}}$: (i), $\tau^{\bar{u}}<T$; (ii), $\tau^{\bar{u}}=T,\ \inf A^{\bar{u}}=T$; (iii) $\tau^{\bar{u}}=T,\ A^{\bar{u}}=\varnothing$. Thus, we need to consider three cases to prove the differentiability of  $\tau^{\bar{u}}$ at $\bar{u}(\cdot)$: (i), $\tau^{\bar{u}}<T$; (ii), $\inf A^{\bar{u}}=T$; (iii), $ A^{\bar{u}}=\varnothing$. Based on the condition $H_0$, for the cases (i), (ii) and  (iii), we can prove that $\tau^{{u}_{\rho}}$ converges to $\tau^{\bar{u}}$ when $\rho\to 0$, respectively. However, we give Example \ref{ex00} to show that Lemma \ref{le-0} and Lemma \ref{le-2} does not right without the condition $H_0$.
}

\begin{lemma}
\label{le-0}
Let Assumptions \ref{ass-b}, \ref{assb-b2} and \ref{ass-fai} hold, and suppose that $h^{\bar{u}}(\tau^{\bar{u}})\neq 0$ and $h^{\bar{u}}(\cdot)$ is continuous at the point $\tau^{\bar{u}}$. We have the following results.

\noindent (i). If $\tau^{\bar{u}}<T$, one obtains
\begin{equation}
\label{le0-e0}
\lim_{\rho\to 0}\left|{\tau^{\bar{u}}-\tau^{u_{\rho}}}\right|
=0.
\end{equation}
(ii). If $\inf\bigg{\{}t:\mathbb{E}[\Phi(\bar{X}(t))]\geq \alpha,\ t\in [0,T] \bigg{\}}=T$, one obtains
\begin{equation}
\lim_{\rho\to 0}\left|{\tau^{\bar{u}}-\tau^{u_{\rho}}}\right|
=0.
\end{equation}
(iii). If $\bigg{\{}t:\mathbb{E}[\Phi(\bar{X}(t))]\geq \alpha,\ t\in [0,T] \bigg{\}}=\varnothing$, we have
\begin{equation}
\lim_{\rho\to 0}\left|{\tau^{\bar{u}}-\tau^{u_{\rho}}}\right|
=0.
\end{equation}

\end{lemma}
\textbf{Proof:} We first  prove {case $i)$}.
Notice that $\tau^{\bar{u}}<T$, we have
\begin{equation}
\label{dech-1}
 \tau^{\bar{u}}=\inf\bigg{\{}t:\mathbb{E}[\Phi(\bar{X}(t))]\geq \alpha,\ t\in [0,T] \bigg{\}}.
\end{equation}
By equation (\ref{mean-0}), it follows that
$$
\mathbb{E}[\Phi(\bar{X}(\tau^{\bar{u}}))]=\Phi(x_0)+\int_0^{\tau^{\bar{u}}}h^{\bar{u}}(t)\mathrm{d}t,\quad
\mathbb{E}[\Phi(X^{u_{\rho}}(\tau^{u_{\rho}}))]=\Phi(x_0)+\int_0^{\tau^{u_{\rho}}}h^{{u_{\rho}}}(t)\mathrm{d}t,
$$
where $u_{\rho}(\cdot)=\bar{u}(\cdot)+\rho v(\cdot)$. For any given $\varepsilon >0$, it follows from Assumptions  (\ref{ass-b}) and (\ref{ass-fai}), there exist $\delta>0$ and $\rho\in (-\delta,\delta)$ such that
\begin{equation}
\label{le1-e1}
\sup_{t\in[0,T]}\left| \mathbb{E}[\Phi(\bar{X}(t)]-\mathbb{E}[\Phi(X^{u_{\rho}}(t)] \right|< \frac{\varepsilon}{2}.
\end{equation}
 By  (\ref{dech-1}) and Remark \ref{re-1}, it follows that
$
\mathbb{E}[\Phi(\bar{X}(\tau^{\bar{u}}))]=\alpha,\ \mathbb{E}[\Phi(\bar{X}(t))]<\alpha,\ t\in [0,\tau^{\bar{u}}).
$
Because $h^{\bar{u}}(\cdot)$ is continuous at the point $\tau^{\bar{u}}$, it follows that
$$
\frac{\text{d}\mathbb{E}[\Phi(\bar{X}(t))]}{\text{d}t}\bigg {|}_{t=\tau^{\bar{u}}}=h^{\bar{u}}(\tau^{\bar{u}}).
$$
Notice that  $h^{\bar{u}}(\tau^{\bar{u}})\neq 0$. Without loss of  generality, we suppose that $h^{\bar{u}}(\tau^{\bar{u}})>0$. Then, there exists $\gamma\in(0,L(\varepsilon))$  such that   $(\tau^{\bar{u}}-\gamma,\tau^{\bar{u}}+\gamma)\subset [0,T]$, where $L(\cdot)>0$ is continuous at $0$, $L(0)=0$, and for $t\in (\tau^{\bar{u}}-\gamma,\tau^{\bar{u}}+\gamma)$, it holds that
$
\left|  \mathbb{E}[\Phi(\bar{X}(t))]-\alpha\right|<\varepsilon,
$
and there exist $t_1\in (\tau^{\bar{u}}-\gamma,\tau^{\bar{u}}),\ t_2\in(\tau^{\bar{u}},\tau^{\bar{u}}+\gamma)$ such that
$$
\alpha-\varepsilon < \mathbb{E}[\Phi(\bar{X}(t_1))]<\alpha-\frac{1}{2}\varepsilon,\quad
\alpha +\frac{1}{2}\varepsilon< \mathbb{E}[\Phi(\bar{X}(t_2))]<\alpha+\varepsilon,
$$
and $\sup_{t\in [0,\tau^{\bar{u}}-\gamma]}\mathbb{E}[\Phi(\bar{X}(t))]<\mathbb{E}[\Phi(\bar{X}(t_1))].$
From equation (\ref{le1-e1}) and $\rho\in (-\delta,\delta)$, we have
\begin{equation*}%
\begin{array}
[l]{ll}%
&\displaystyle \sup_{t\in [0,\tau^{\bar{u}}-\gamma]}\mathbb{E}[\Phi(X^{{u}_\rho}(t))]\\
\leq&  \displaystyle\sup_{t\in [0,\tau^{\bar{u}}-\gamma]}\mathbb{E}[\Phi(\bar{X}(t))]+\frac{\varepsilon}{2}\\
< &\displaystyle \mathbb{E}[\Phi(\bar{X}(t_1))]+\frac{\varepsilon}{2}\\
< &\alpha,
\end{array}
\end{equation*}%
and
$$
\alpha<\mathbb{E}[\Phi(\bar{X}(t_2))]-\frac{\varepsilon}{2} \leq \mathbb{E}[\Phi(X^{u_{\rho}}(t_2))].
$$
This implies that
$$
\tau^{u_{\rho}}\in (\tau^{\bar{u}}-\gamma,t_2)\subset (\tau^{\bar{u}}-\gamma,\tau^{\bar{u}}+\gamma),
$$
and thus $\left|\tau^{u_{\rho}}-\tau^{\bar{u}}\right| <L(\varepsilon)$.

{ Second, we consider the case $(ii)$. Notice that $\inf\bigg{\{}t:\mathbb{E}[\Phi(\bar{X}(t))]\geq \alpha,\ t\in [0,T] \bigg{\}}=T$, for sufficiently small $\rho$, we have $\tau^{{u_{\rho}}}\leq T$. Similar with the proof of case $(i)$, we can obtain $\displaystyle \lim_{\rho\to 0}\left|{\tau^{\bar{u}}-\tau^{u_{\rho}}}\right|=0.$
In the end, we consider case $(iii)$. Notice that $\bigg{\{}t:\mathbb{E}[\Phi(\bar{X}(t))]\geq \alpha,\ t\in [0,T] \bigg{\}}=\varnothing$, thus for sufficiently small $\rho$, we have $\tau^{{u_{\rho}}}=T$ and
$
\lim_{\rho\to 0}\left|{\tau^{\bar{u}}-\tau^{u_{\rho}}}\right|=0.
$
}This completes the proof. $\ \ \ \ \ \ \ \ \Box$
\bigskip

Let $y^{}(\cdot)$  be the solution
of the following variational equation:%
\begin{equation}
\label{apro-1}
\begin{array}
[c]{rl}%
\text{d}{y}(t)= & \big{[}b_{x}(\bar{X}{(t)},\bar{u}(t))y(t)+b_{u}(\bar{X}{(t)},\bar{u}(t)) v(t)\big{]}\text{d}t\\
&+ \displaystyle \sum_{j=1}^d\big{[} \sigma_{x}^j(\bar{X}{(t)},\bar{u}(t))y(t)+\sigma^j_u(\bar{X}{(t)},\bar{u}(t)) v(t)\big{]}\text{d}W^j(t), \\
y(0)= & 0,\quad t\in (0,T].
\end{array}
\end{equation}
{ The following lemma is classical, and we omit the proof, see \cite{B81} and \cite{B78}.}
\begin{lemma}
\label{le-1} Let Assumptions \ref{ass-b}, \ref{assb-b2} and \ref{ass-fai} hold. We have
\begin{equation}%
\begin{array}
[l]{l}%
\displaystyle\lim_{\rho\to 0}\displaystyle\sup_{t\in \lbrack0,T]}\mathbb{E}\left|\rho^{-1} (X^{\rho}(t)-\bar{X}(t))-y(t)\right| =0.  \\
\end{array}
 \label{var-1}%
\end{equation}
\end{lemma}
\bigskip

\begin{remark}
We consider a special case that $\Phi(x)=x,\ m=1$. Thus
\begin{equation}
\mathbb{E}[X^u(s)]=x_0+\int_0^sh^u(t)\mathrm{d}t,
\end{equation}
where $h^{{u}}(t)=\mathbb{E}\big{[}b(X^{{u}}(t),{u}(t))\big{]}$.
It follows from equation (\ref{var-1}) that
\begin{equation}
\label{re-e1}
\begin{array}
[c]{rl}%
&\displaystyle \lim_{\rho\to 0}\frac{h^{u_{\rho}}(t)-h^{\bar{u}}(t)}{\rho} =\mathbb{E}\bigg{[}b_{x}(\bar{X}{(t)},\bar{u}(t))y(t)
+b_{u}(\bar{X}{(t)},\bar{u}(t))v(t)\bigg{]}.
\end{array}
\end{equation}
For notation simplicity, we set $ \bar{h}(v(t),t)=\displaystyle \lim_{\rho\to 0}\frac{h^{u_{\rho}}(t)-h^{\bar{u}}(t)}{\rho}$. For general $\Phi(\cdot)$, we give the explicit formula for $\bar{h}(v(t),t)$ in Appendix \ref{app-a}.
\end{remark}
\begin{lemma}
\label{le-2}
Let Assumptions \ref{ass-b}, \ref{assb-b2} and \ref{ass-fai} hold. Suppose that $h^{\bar{u}}(\tau^{\bar{u}})\neq 0$, and $h^{\bar{u}}(\cdot)$ is continuous at the point $\tau^{\bar{u}}$. We have the following results.

\noindent (i). If $\tau^{\bar{u}}<T$, one obtains
\begin{equation}
\label{le2-e0}
\lim_{\rho\to 0}\frac{\tau^{\bar{u}}-\tau^{u_{\rho}}}{\rho}
=\int_0^{\tau^{\bar{u}}}\frac{ \bar{h}(v(t),t)}{h^{\bar{u}}(\tau^{\bar{u}})}\mathrm{d}t.
\end{equation}
(ii). If $\inf\bigg{\{}t:\mathbb{E}[\Phi(\bar{X}(t))]\geq \alpha,\ t\in [0,T] \bigg{\}}=T$, then there exists sequence $\rho_n\to 0$ as $n\to +\infty$ such that
\begin{equation}
\lim_{n\to +\infty}\frac{\tau^{\bar{u}}-\tau^{u_{\rho_n}}}{\rho_n}=
\int_0^{\tau^{\bar{u}}}\frac{ \bar{h}(v(t),t)}{h^{\bar{u}}(\tau^{\bar{u}})}\mathrm{d}t\ \  \mathrm{or}\ \  0.
\end{equation}
(iii). If $\bigg{\{}t:\mathbb{E}[\Phi(\bar{X}(t))]\geq \alpha,\ t\in [0,T] \bigg{\}}=\varnothing$, we have
\begin{equation}
\lim_{\rho\to 0}\frac{\tau^{\bar{u}}-\tau^{u_{\rho}}}{\rho}=0.
\end{equation}

\end{lemma}
\textbf{Proof:} We first prove case $(i)$. Notice that for $ \tau^{\bar{u}}<T$,
\begin{equation*}
 \tau^{\bar{u}}=\inf\bigg{\{}t:\mathbb{E}[\Phi({\bar{X}}(t))]\geq \alpha,\ t\in [0,T] \bigg{\}}.
\end{equation*}
By $(i)$ of Lemma \ref{le-0}, we have that $\displaystyle\lim_{\rho\to 0}\left|\tau^{u_{\rho}}-\tau^{\bar{u}}\right|=0$, for small sufficiently  $\rho$, it follows that
\begin{equation*}
 \tau^{{u_{\rho}}}=\inf\bigg{\{}t:\mathbb{E}[\Phi(X^{{u_{\rho}}}(t))]\geq \alpha,\ t\in [0,T] \bigg{\}}.
\end{equation*}
This implies that
$
\mathbb{E}[\Phi({\bar{X}}(\tau^{\bar{u}}))]=\mathbb{E}[\Phi(X^{{u_{\rho}}}(\tau^{u_{\rho}}))]=\alpha.
$
Combining
$$
\mathbb{E}[\Phi(\bar{X}(\tau^{\bar{u}}))]=\Phi(x_0)+
\int_0^{\tau^{\bar{u}}}h^{\bar{u}}(t)\mathrm{d}t,
$$
and
$$
\mathbb{E}[\Phi(X^{u_{\rho}}(\tau^{u_{\rho}}))]=\Phi(x_0)
+\int_0^{\tau^{u_{\rho}}}h^{{u_{\rho}}}(t)\mathrm{d}t,
$$
we have
$$
\int_0^{\tau^{\bar{u}}}h^{\bar{u}}(t)\mathrm{d}t=\int_0^{\tau^{u_{\rho}}}h^{{u_{\rho}}}(t)\mathrm{d}t,
$$
and
$$
\int_{\tau^{u_{\rho}}}^{\tau^{\bar{u}}}h^{\bar{u}}(t)\text{d}t
=\int_0^{\tau^{u_{\rho}}}\bigg{[}h^{{u_{\rho}}}(t)-h^{\bar{u}}(t)\bigg{]}\text{d}t.
$$
Dividing  on both sides of the above equation by $\rho$, we obtain
$$
\displaystyle \lim_{\rho\to 0}\frac{\int_{\tau^{u_{\rho}}}^{\tau^{\bar{u}}}h^{\bar{u}}(t)\text{d}t}{\rho}
=\displaystyle \lim_{\rho\to 0}\int_0^{\tau^{u_{\rho}}}\frac{h^{{u_{\rho}}}(t)-h^{\bar{u}}(t)}{\rho}\text{d}t.
$$

Again by $(i)$ of Lemma \ref{le-0}, $\displaystyle\lim_{\rho\to 0}\left|\tau^{u_{\rho}}-\tau^{\bar{u}}\right|=0$, we have
\begin{equation}
\label{le2-e1}
\begin{array}
[c]{rl}%
& \displaystyle \lim_{\rho\to 0}\frac{\int_{\tau^{u_{\rho}}}^{\tau^{\bar{u}}}h^{\bar{u}}(t)\text{d}t}{\rho}
= \displaystyle \lim_{\rho\to 0}\frac{\tau^{\bar{u}}-\tau^{u_{\rho}}}{\rho}\bigg{[}h^{\bar{u}}(\tau^{\bar{u}})
+o(1)\bigg{]},\\
\end{array}
\end{equation}
where $o(1)$ converges to $0$ as $\left|\tau^{u_{\rho}}-\tau^{\bar{u}} \right|\to 0$. Note that
$$
\bar{h}(v(t),t)=\lim_{\rho\to 0}\frac{h^{u_{\rho}}(t)-h^{\bar{u}}(t)}{\rho},\ t\in[0,T].
$$
It follows that
\begin{equation}
\label{le2-e2}
\begin{array}
[c]{rl}%
& \displaystyle \lim_{\rho\to 0}\int_0^{\tau^{u_{\rho}}} \frac{h^{u_{\rho}}(t)-h^{\bar{u}}(t)}{\rho}\text{d}t
=\displaystyle \int_0^{\tau^{\bar{u}}} \bar{h}(v(t),t)  \text{d}t.
\end{array}
\end{equation}
Combining equations (\ref{le2-e1}) and (\ref{le2-e2}), we obtain
$$
\lim_{\rho\to 0}\frac{\tau^{\bar{u}}-\tau^{u_{\rho}}}{\rho}
=\int_0^{\tau^{\bar{u}}}\frac{ \bar{h}(v(t),t)}{h^{\bar{u}}(\tau^{\bar{u}})}\text{d}t.
$$

{ Second, we consider the case $(ii)$. Notice that $\inf\bigg{\{}t:\mathbb{E}[\Phi(\bar{X}(t))]\geq \alpha,\ t\in [0,T] \bigg{\}}=T$, if there exists sequence $\rho_n\to 0$ as $n\to +\infty$ such that
$$
\tau^{{u_{\rho_n}}}=\inf\bigg{\{}t:\mathbb{E}[\Phi(X^{{u_{\rho_n}}}(t))]\geq \alpha,\ t\in [0,T] \bigg{\}}<T.
$$
Similar with the proof of case $(i)$, by the case $(ii)$ of Lemma \ref{le-0}, we can obtain
\begin{equation*}
\lim_{n\to +\infty}\frac{\tau^{\bar{u}}-\tau^{u_{\rho_n}}}{\rho_n}=
\int_0^{\tau^{\bar{u}}}\frac{ \bar{h}(v(t),t)}{h^{\bar{u}}(\tau^{\bar{u}})}\mathrm{d}t.
\end{equation*}
If there exists sequence $\rho_n\to 0$ as $n\to +\infty$ such that
$$
\inf\bigg{\{}t:\mathbb{E}[\Phi(X^{{u_{\rho_n}}}(t))]\geq \alpha,\ t\in [0,T] \bigg{\}}=+\infty,
$$
then, $\tau^{{u_{\rho_n}}}=T$, and
\begin{equation*}
\lim_{n\to +\infty}\frac{\tau^{\bar{u}}-\tau^{u_{\rho_n}}}{\rho_n}=0.
\end{equation*}
Thus, the case $(ii)$ is right.

In the end, we consider case $(iii)$. Notice that $\bigg{\{}t:\mathbb{E}[\Phi(\bar{X}(t))]\geq \alpha,\ t\in [0,T] \bigg{\}}=\varnothing$, thus for sufficiently small $\rho$,
$$
\inf\bigg{\{}t:\mathbb{E}[\Phi(X^{{u_{\rho}}}(t))]\geq \alpha,\ t\in [0,T] \bigg{\}}=+\infty,
$$
and $\tau^{{u_{\rho}}}=T$, which shows that case $(iii)$ is right.
}
\noindent This completes the proof. $\ \ \ \ \ \ \ \ \Box$

\bigskip
{ In the following example, we show that why we need the conditions that $h^{\bar{u}}(\tau^{\bar{u}})\neq 0$ and $h^{\bar{u}}(\cdot)$ is continuous at the point $\tau^{\bar{u}}$ in Lemma \ref{le-2}.}
\begin{example}\label{ex00}
Let { $T=2,\ U=[-2,2]$}, $m=d=1$, $b(x,u)=u,\ \sigma(x,u)=0$ and $\Phi(x)=x$.

 \textbf{Case 1):} For a given cost functional,
$$
J(u(\cdot))=\int_0^1(u(t)-1)^2\mathrm{d}t+\int_1^2(u(t)-0.5)^2\mathrm{d}t,
$$
an optimal pair is given as follows:
\begin{equation}
   (\bar{u}(t),\bar{X}(t))=
   \begin{cases}
    (1,t), &\mbox{ $0\leq t\leq 1$,}\\
  (0.5,0.5+0.5t), &\mbox{ $ 1< t\leq 2$. }
   \end{cases}
  \end{equation}
The shape  of $\bar{X}(\cdot)$ is illustrated in Figure 2.
\begin{figure}[h]
\begin{center}
\includegraphics[width=4.8 in]{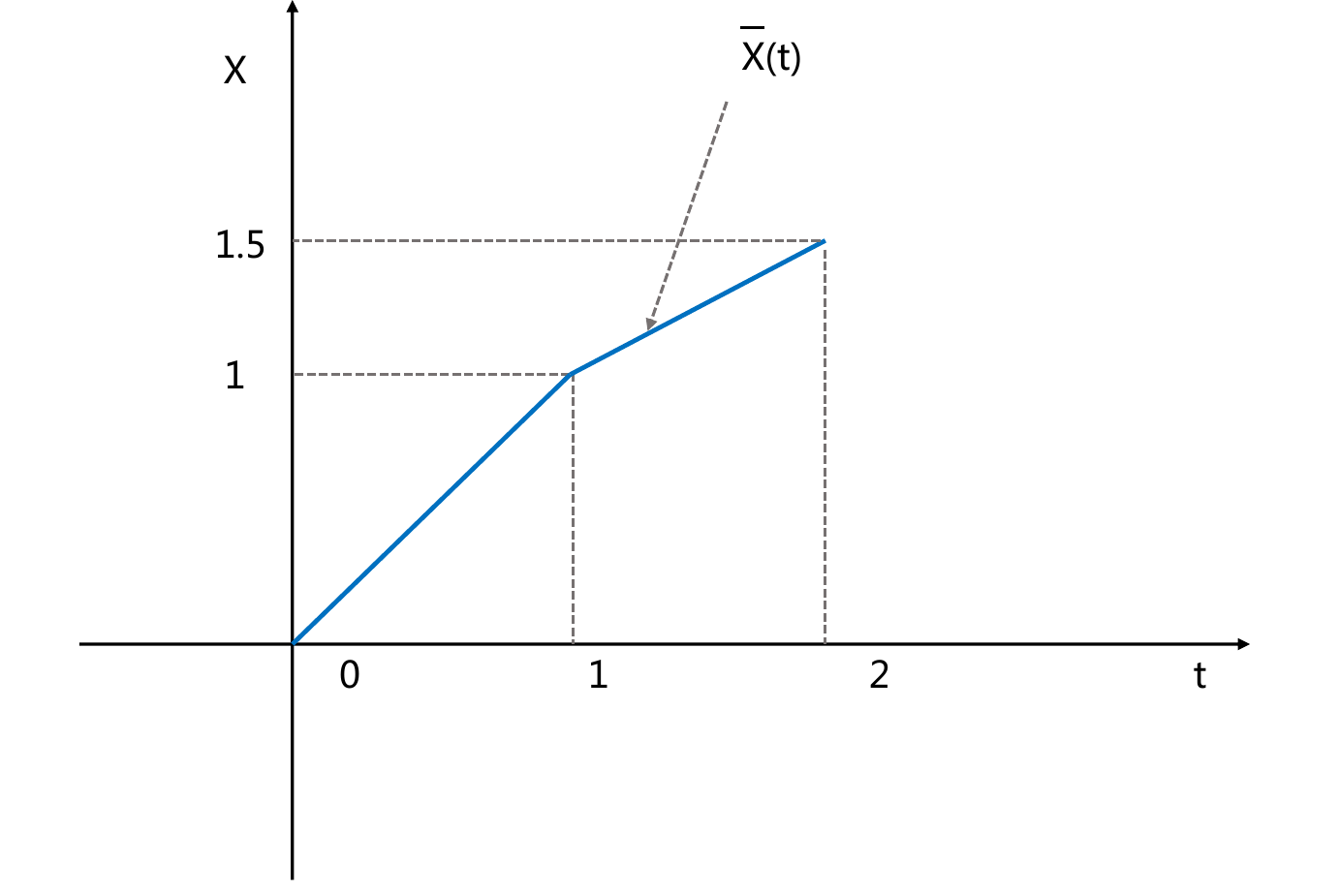}%
\caption{The shape  of $\bar{X}(\cdot)$}
\end{center}
 \label{fig11}
\end{figure}

We take $\alpha=1$. Thus, $\tau^{\bar{u}}=1$. Let $v=1$. Then,  for $t\in [0,2]$ the control $u_{\rho}(\cdot)$ is given as follows:
$$
u_{\rho}(t)= \bar{u}(t)+\rho.
$$
We calculate $\tau^{u_{\rho}}$ as
\begin{equation}
  \tau^{u_{\rho}}=
   \begin{cases}
\displaystyle  \frac{1}{1+\rho}, &\mbox{ $\rho>0$,}\\
\displaystyle   \frac{0.5}{0.5+\rho}, &\mbox{ $ \rho<0$,}
   \end{cases}
  \end{equation}
and
$$
\lim_{\rho\to 0^{+}}\frac{\tau^{\bar{u}}-\tau^{u_{\rho}}}{\rho}=\lim_{\rho\to 0^{+}}\frac{1-\frac{1}{1+\rho}}{\rho}=1,
$$
$$
\lim_{\rho\to 0^{-}}\frac{\tau^{\bar{u}}-\tau^{u_{\rho}}}{\rho}=\lim_{\rho\to 0^{-}}\frac{1-\frac{0.5}{0.5+\rho}}{\rho}=2.
$$
This shows that $\displaystyle \lim_{\rho\to 0}\frac{\tau^{\bar{u}}-\tau^{u_{\rho}}}{\rho}$ does not exist. This is because  $h^{\bar{u}}(\cdot)=\bar{u}(\cdot)$ is not continuous at the point $\tau^{\bar{u}}=1$. This example shows that we cannot deal with a form of $\bar{X}(\cdot)$ that is similar to that in Figure 2.

\textbf{Case 2):} For a given cost functional,
$$
J(u(\cdot))=\int_0^2(u(t)-2+2t)^2\mathrm{d}t,
$$
an optimal pair is given as follows:
\begin{equation}
   (\bar{u}(t),\bar{X}(t))=
    (2-2t,2t-t^2), \ 0\leq t\leq 2.
  \end{equation}
The shape of $\bar{X}(\cdot)$ is illustrated in Figure 3.
\begin{figure}[h] \label{fig2}
\begin{center}
\includegraphics[width=4.0 in]{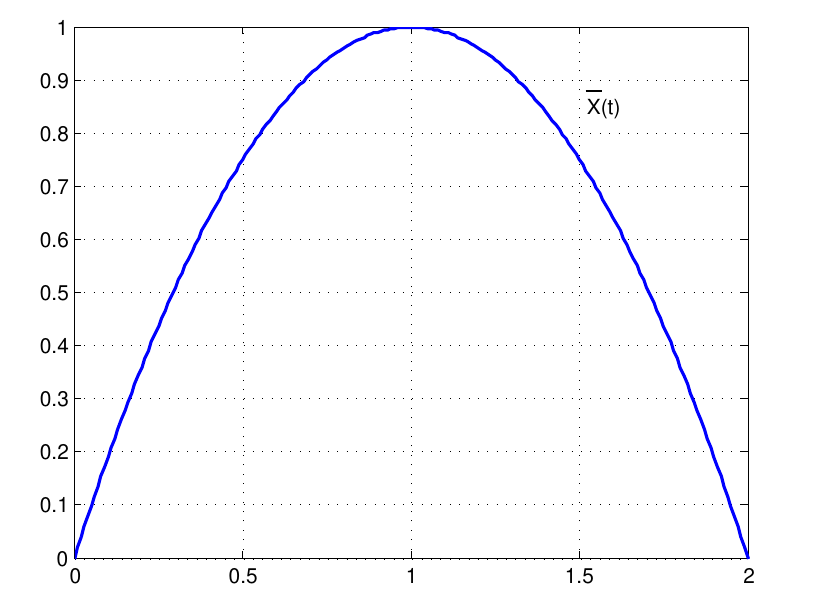}
 \caption{The shape  of $\bar{X}(\cdot)$  }
\end{center}
\end{figure}

We take $\alpha=1$, and thus $\tau^{\bar{u}}=1$. Let $v=1$. Then, for $t\in [0,2]$, the control $u_{\rho}(\cdot)$ is given as follows:
$$
u_{\rho}(t)= \bar{u}(t)+\rho.
$$
A simple calculation shows that $\tau^{u_{\rho}}$ is given by
\begin{equation}
\displaystyle    \tau^{u_{\rho}}=
   \begin{cases}
\frac{2+\rho-\sqrt{4\rho+\rho^2}}{2}, &\mbox{ $\rho>0$,}\\
 2, &\mbox{ $ \rho<0$,}
   \end{cases}
  \end{equation}
and
\begin{equation*}
\begin{array}
[c]{rl}
&\displaystyle \lim_{\rho\to 0^{+}}\frac{\tau^{\bar{u}}-\tau^{u_{\rho}}}{\rho}=\lim_{\rho\to 0^{+}}\frac{1-\frac{2+\rho-\sqrt{4\rho+\rho^2}}{2}}{\rho}=+\infty,\\
&\displaystyle \lim_{\rho\to 0^{-}}\frac{\tau^{\bar{u}}-\tau^{u_{\rho}}}{\rho}=\lim_{\rho\to 0^{+}}\frac{1-2}{\rho}=-\infty.
\end{array}
\end{equation*}
Thus, $\displaystyle\lim_{\rho\to 0}\frac{\tau^{\bar{u}}-\tau^{u_{\rho}}}{\rho}$ does not exist. This is because  $h^{\bar{u}}(\tau^{\bar{u}})=2-2\tau^{\bar{u}}=0$. This example shows that we cannot deal with a form of $\bar{X}(\cdot)$ that is similar to that in Figure 3.
\end{example}

We now derive the variational equation for cost functional (\ref{cost-1}) in the following lemma.
\begin{lemma}
\label{le-3} Let Assumptions \ref{ass-b}, \ref{assb-b2} and \ref{ass-fai} hold and suppose that $h^{\bar{u}}(\tau^{\bar{u}})\neq 0$ and $h^{\bar{u}}(\cdot)$ is continuous at the point $\tau^{\bar{u}}$. We have the following results.

\noindent (i). If $\tau^{\bar{u}}<T$, one obtains
\begin{equation}
\label{le3-e0}
\begin{array}
[c]{rl}
& \rho^{-1}\left[J(u_{\rho}(\cdot))-J(\bar{u}(\cdot))\right] \\
= &\displaystyle -\int_0^{\tau^{\bar{u}}}\bigg{[}\frac{ \tilde{\Psi}^{\bar{u}}(\tau^{\bar{u}})\bar{h}(v(t),t)}{h^{\bar{u}}(\tau^{\bar{u}})}+\frac{ \mathbb{E}[f(\bar{X}(\tau^{\bar{u}}),\bar{u}(\tau^{\bar{u}}))]\bar{h}(v(t),t)}{h^{\bar{u}}(\tau^{\bar{u}})}\bigg{]}\mathrm{d}t
+\mathbb{E}\bigg{[} \Psi_x(\bar{X}(\tau^{\bar{u}}))^{\top}y(\tau^{\bar{u}})\big{)}\bigg{]}\\
&+\displaystyle \mathbb{E}\int_{0}^{\tau^{\bar{u}}}\bigg{[}
f_x(\bar{X}{(t)},\bar{u}(t))^{\top}y(t)+f_u(\bar{X}{(t)},\bar{u}(t))^{\top}v(t)  \bigg{]}\mathrm{d}t+o(1),
\end{array}
\end{equation}
where
$$
\tilde{\Psi}^{\bar{u}}(\tau^{\bar{u}})=\mathbb{E}\bigg{[}\Psi_x(\bar{X}(\tau^{\bar{u}}))^{\top}b(\bar{X}(\tau^{\bar{u}}),\bar{u}(\tau^{\bar{u}}))
+\frac{1}{2}\sum_{j=1}^d\sigma^j(\bar{X}(\tau^{\bar{u}}),\bar{u}(\tau^{\bar{u}}))^{\top}\Psi_{xx}(\bar{X}(\tau^{\bar{u}}))
\sigma^j(\bar{X}(\tau^{\bar{u}}),\bar{u}(\tau^{\bar{u}}))\bigg{]}.
$$
(ii). If $\inf\bigg{\{}t:\mathbb{E}[\Phi(\bar{X}(t))]\geq \alpha,\ t\in [0,T] \bigg{\}}=T$, one obtains
\begin{equation}
\begin{array}
[c]{rl}
& \rho^{-1}\left[J(u_{\rho}(\cdot))-J(\bar{u}(\cdot))\right] \\
= &\displaystyle -\int_0^{\tau^{\bar{u}}}\bigg{[}\frac{ \tilde{\Psi}^{\bar{u}}(\tau^{\bar{u}})\bar{h}(v(t),t)}{h^{\bar{u}}(\tau^{\bar{u}})}+\frac{ \mathbb{E}[f(\bar{X}(\tau^{\bar{u}}),\bar{u}(\tau^{\bar{u}}))]\bar{h}(v(t),t)}{h^{\bar{u}}(\tau^{\bar{u}})}\bigg{]}\mathrm{d}t
+\mathbb{E}\bigg{[} \Psi_x(\bar{X}(\tau^{\bar{u}}))^{\top}y(\tau^{\bar{u}})\big{)}\bigg{]}\\
&+\displaystyle\int_{0}^{\tau^{\bar{u}}}\mathbb{E}\bigg{[}
f_x(\bar{X}{(t)},\bar{u}(t))^{\top}y(t)+f_u(\bar{X}{(t)},\bar{u}(t))^{\top}v(t)  \bigg{]}\mathrm{d}t+o(1),
\end{array}
\end{equation}
or
\begin{equation}
\begin{array}
[c]{rl}
& \rho^{-1}\left[J(u_{\rho}(\cdot))-J(\bar{u}(\cdot))\right] \\
= &\displaystyle \mathbb{E}\bigg{[} \Psi_x(\bar{X}(\tau^{\bar{u}}))^{\top}y(\tau^{\bar{u}})\big{)}\bigg{]}+\displaystyle\int_{0}^{\tau^{\bar{u}}}\mathbb{E}\bigg{[}
f_x(\bar{X}{(t)},\bar{u}(t))^{\top}y(t)+f_u(\bar{X}{(t)},\bar{u}(t))^{\top}v(t)  \bigg{]}\mathrm{d}t+o(1).
\end{array}
\end{equation}
(iii). If $\bigg{\{}t:\mathbb{E}[\Phi(\bar{X}(t))]\geq \alpha,\ t\in [0,T] \bigg{\}}=\varnothing$, we have
\begin{equation}
\begin{array}
[c]{rl}
& \rho^{-1}\left[J(u_{\rho}(\cdot))-J(\bar{u}(\cdot))\right] \\
= &\mathbb{E}\bigg{[} \Psi_x(\bar{X}(\tau^{\bar{u}}))^{\top}y(\tau^{\bar{u}})\big{)}\bigg{]}
+\displaystyle\int_{0}^{\tau^{\bar{u}}}\mathbb{E}\bigg{[}
f_x(\bar{X}{(t)},\bar{u}(t))^{\top}y(t)+f_u(\bar{X}{(t)},\bar{u}(t))^{\top}v(t)  \bigg{]}\mathrm{d}t+o(1).
\end{array}
\end{equation}
\end{lemma}
\textbf{Proof: }
We first prove case $(i)$. Notice that%
\begin{equation}%
\begin{array}
[c]{rl}
& J(u_{\rho}(\cdot))-J(\bar{u}(\cdot))\\
= &\mathbb{E}\bigg{[} \Psi(X^{\rho}(\tau^{u_{\rho}}))-\Psi(\bar{X}(\tau^{\bar{u}}))+{\displaystyle \int_{0}^{\tau^{u_{\rho}}}}
f(X{}^{\rho}(t),u_{\rho}(t))\text{d}t-{\displaystyle \int_{0}^{\tau^{\bar{u}}}}f(\bar{X}{(t)},\bar{u}(t))\text{d}t\bigg{]},\\
=&\mathbb{E}\bigg{[} \Psi(X^{\rho}(\tau^{u_{\rho}}))-\Psi(X^{\rho}(\tau^{\bar{u}}))\bigg{]}
+\mathbb{E}\bigg{[} \displaystyle \int_{\tau^{\bar{u}}}^{\tau^{u_{\rho}}}
f(X^{\rho}(t),u_{\rho}(t))\text{d}t\bigg{]}\\
&+\mathbb{E}\bigg{[} \Psi(X^{\rho}(\tau^{\bar{u}}))-\Psi(\bar{X}(\tau^{\bar{u}}))\bigg{]}+{\displaystyle \mathbb{E}\int_{0}^{\tau^{\bar{u}}}}\bigg{[}f(X^{\rho}(t),u_{\rho}(t))-
f(\bar{X}{(t)},\bar{u}(t))\bigg{]}\text{d}t.
\end{array}
\label{value-0}%
\end{equation}
Denoting
\begin{equation*}%
\begin{array}
[c]{rl}
&I_1:=\rho^{-1}\mathbb{E}\bigg{[} \Psi(X^{\rho}(\tau^{u_{\rho}}))-\Psi(X^{\rho}(\tau^{\bar{u}}))\bigg{]},\\
&I_2:=\rho^{-1}\mathbb{E}\bigg{[} \displaystyle \int_{\tau^{\bar{u}}}^{\tau^{u_{\rho}}}
f(X^{\rho}(t),u_{\rho}(t))\text{d}t\bigg{]},\\
&I_3:=\rho^{-1}\mathbb{E}\bigg{[} \Psi(X^{\rho}(\tau^{\bar{u}}))-\Psi(\bar{X}(\tau^{\bar{u}}))\bigg{]},\\
&I_4:=\rho^{-1}{\displaystyle\mathbb{E} \int_{0}^{\tau^{\bar{u}}}}\bigg{[}f(X^{\rho}(t),u_{\rho}(t))-
f(\bar{X}{(t)},\bar{u}(t))\bigg{]}\text{d}t,
\end{array}
\end{equation*}
we first consider the terms $I_1$ and $I_2$. Applying It\^{o} formula to $\Psi(\cdot)$, we have
\begin{equation}%
\begin{array}
[c]{ll}%
& \text{d}\Psi(X^{\rho}(t))=\Psi_x(X^{\rho}(t))^{\top}\text{d}X^{\rho}(t)
+\displaystyle \frac{1}{2}\big{[}\text{d}X^{\rho}(t)\big{]}^{\top}\Psi_{xx}(X^{\rho}(t))\text{d}X^{\rho}(t).\\
\end{array}
\end{equation}
By equation (\ref{ODE_1}), it follows that
\begin{equation*}%
\begin{array}
[c]{ll}%
 \text{d}\Psi(X^{\rho}(t))=&\bigg{[}\Psi_x(X^{\rho}(t))^{\top}b(X^{\rho}(t),u_{\rho}(t))
 +\displaystyle \frac{1}{2}\sum_{j=1}^d\sigma^j(X^{\rho}(t),u_{\rho}(t))^{\top}
 \Psi_{xx}(X^{\rho}(t))\sigma^j(X^{\rho}(t),u_{\rho}(t))\bigg{]}\text{d}t\\
&+\Psi_x(X^{\rho}(t))^{\top}\sigma(X^{\rho}(t),u_{\rho}(t))\text{d}W(t),\\
\end{array}
\end{equation*}
and thus
\begin{equation*}%
\begin{array}
[c]{ll}%
 &\Psi(X^{\rho}(\tau^{u_{\rho}}))-\Psi(X^{\rho}(\tau^{\bar{u}}))\\
=&\displaystyle \int_{\tau^{\bar{u}}}^{\tau^{u_{\rho}}}\bigg{[}\Psi_x(X^{\rho}(t))^{\top}b(X^{\rho}(t),u_{\rho}(t))
+\frac{1}{2}\sum_{j=1}^d\sigma^j(X^{\rho}(t),u_{\rho}(t))^{\top}\Psi_{xx}(X^{\rho}(t))
\sigma^j(X^{\rho}(t),u_{\rho}(t))\bigg{]}\text{d}t\\
&+\displaystyle\int_{\tau^{\bar{u}}}^{\tau^{u_{\rho}}}\Psi_x(X^{\rho}(t))^{\top}\sigma(X^{\rho}(t),u_{\rho}(t))\text{d}W(t),\\
\end{array}
\end{equation*}
and
\begin{equation*}%
\begin{array}
[c]{ll}%
 &\mathbb{E}\big{[}\Psi(X^{\rho}(\tau^{u_{\rho}}))-\Psi(X^{\rho}(\tau^{\bar{u}}))\big{]}\\
=&\displaystyle \int_{\tau^{\bar{u}}}^{\tau^{u_{\rho}}}\mathbb{E}\bigg{[}\Psi_x(X^{\rho}(t))^{\top}b(X^{\rho}(t),u_{\rho}(t))
+\frac{1}{2}\sum_{j=1}^d\sigma^j(X^{\rho}(t),u_{\rho}(t))^{\top}\Psi_{xx}(X^{\rho}(t))
\sigma^j(X^{\rho}(t),u_{\rho}(t))\bigg{]}\text{d}t.\\
\end{array}
\end{equation*}
For $t\in [0,T]$, we set
\begin{equation*}%
\begin{array}
[c]{ll}%
&\displaystyle\tilde{\Psi}^{u_{\rho}}(t)=\mathbb{E}\bigg{[}\Psi_x(X^{\rho}(t))^{\top}b(X^{\rho}(t),u_{\rho}(t))
+\frac{1}{2}\sum_{j=1}^d\sigma^j(X^{\rho}(t),u_{\rho}(t))^{\top}\Psi_{xx}(X^{\rho}(t))
\sigma^j(X^{\rho}(t),u_{\rho}(t))\bigg{]},
\end{array}
\end{equation*}
and thus
\begin{equation}
\label{le3-mean}
\mathbb{E}\big{[}\Psi(X^{\rho}(\tau^{u_{\rho}}))-
\Psi(X^{\rho}(\tau^{\bar{u}}))\big{]}=\int_{\tau^{\bar{u}}}^{\tau^{u_{\rho}}} \tilde{\Psi}^{u_{\rho}}(t)\text{d}t.
\end{equation}
Applying $(i)$ of Lemma \ref{le-2}, it follows from Assumption \ref{ass-fai} and equation (\ref{le2-e0}) that
\begin{equation}
\label{le3-e1}
I_1={\rho^{-1}}\mathbb{E}\big{[}\Psi(X^{\rho}(\tau^{u_{\rho}}))-\Psi(X^{\rho}(\tau^{\bar{u}}))\big{]}
=\int_0^{\tau^{\bar{u}}}-\frac{ \tilde{\Psi}^{\bar{u}}(\tau^{\bar{u}})\bar{h}(v(t),t)}{h^{\bar{u}}(\tau^{\bar{u}})}\text{d}t+o(1).
\end{equation}
Similarly, we obtain
\begin{equation}
\label{le3-e2}
I_2=\rho^{-1}\mathbb{E}\bigg{[} \displaystyle \int_{\tau^{\bar{u}}}^{\tau^{u_{\rho}}}
f(X^{\rho}(t),u_{\rho}(t))\text{d}t\bigg{]}=\int_0^{\tau^{\bar{u}}}-\frac{ \mathbb{E}[f(\bar{X}(\tau^{\bar{u}}),\bar{u}(\tau^{\bar{u}}))]\bar{h}(v(t),t)}{h^{\bar{u}}(\tau^{\bar{u}})}\text{d}t+o(1).
\end{equation}

In the following, we consider $I_3$ and $I_4$:
\begin{equation*}%
\begin{array}
[c]{ll}%
I_3&=\rho^{-1}\mathbb{E}\bigg{[} \Psi(X^{\rho}(\tau^{\bar{u}}))-\Psi(\bar{X}(\tau^{\bar{u}}))\bigg{]}\\
&=\rho^{-1}\mathbb{E}\bigg{[} \Psi_x(X^{\rho}(\tau^{\bar{u}}))^{\top}\big{(}X^{\rho}
(\tau^{\bar{u}})-\bar{X}(\tau^{\bar{u}})\big{)}\bigg{]}+o(1)\\
&=\mathbb{E}\bigg{[} \Psi_x(X^{\rho}(\tau^{\bar{u}}))^{\top}\big{[}\displaystyle \frac{X^{\rho}
(\tau^{\bar{u}})-\bar{X}(\tau^{\bar{u}})}{\rho}-y(\tau^{\bar{u}})+y(\tau^{\bar{u}})\big{]}\bigg{]}+o(1).\\
\end{array}
\end{equation*}
Applying Lemma \ref{le-1}, it follows from  (\ref{var-1}) that
\begin{equation}
\label{le3-e3}
I_3=\mathbb{E}\bigg{[} \Psi_x(\bar{X}(\tau^{\bar{u}}))^{\top}y(\tau^{\bar{u}})\bigg{]}+o(1).
\end{equation}
By equation (\ref{var-1}) and $u_{\rho}(\cdot)=\bar{u}(\cdot)+\rho v(\cdot)$, we have
\begin{equation*}%
\begin{array}
[c]{ll}%
I_4&=\rho^{-1}{\displaystyle \int_{0}^{\tau^{\bar{u}}}}\mathbb{E}\bigg{[}f(X^{\rho}(t),u_{\rho}(t))-
f(\bar{X}{(t)},\bar{u}(t))\bigg{]}\text{d}t\\
&=\rho^{-1}\displaystyle \int_{0}^{\tau^{\bar{u}}}\mathbb{E}\bigg{[}f_x(\bar{X}{(t)},\bar{u}(t))^{\top}\big{(}X^{\rho}
(t)-\bar{X}(t)\big{)}+f_u(\bar{X}{(t)},\bar{u}(t))^{\top}\big{(}u_{\rho}
(t)-\bar{u}(t)\big{)}\bigg{]}\text{d}t+o(1)\\
&=\displaystyle \int_{0}^{\tau^{\bar{u}}}\mathbb{E}\bigg{[}
f_x(\bar{X}{(t)},\bar{u}(t))^{\top}y(t)+f_u(\bar{X}{(t)},\bar{u}(t))^{\top}v(t)  \bigg{]}\text{d}t+o(1),\\
\end{array}
\end{equation*}
which implies that
\begin{equation}
\label{le3-e4}
I_4=\displaystyle \int_{0}^{\tau^{\bar{u}}}\mathbb{E}\bigg{[}
f_x(\bar{X}{(t)},\bar{u}(t))^{\top}y(t)+f_u(\bar{X}{(t)},\bar{u}(t))^{\top}v(t)  \bigg{]}\text{d}t+o(1).
\end{equation}
Combining equations (\ref{le3-e1}), (\ref{le3-e2}), (\ref{le3-e3}) and (\ref{le3-e4}), we obtain
\begin{equation*}
\begin{array}
[c]{rl}
& \rho^{-1}\left[J(u_{\rho}(\cdot))-J(\bar{u}(\cdot))\right] \\
= &\displaystyle -\int_0^{\tau^{\bar{u}}}\bigg{[}\frac{ \tilde{\Psi}^{\bar{u}}(\tau^{\bar{u}})\bar{h}(v(t),t)}{h^{\bar{u}}(\tau^{\bar{u}})}+\frac{ \mathbb{E}[f(\bar{X}(\tau^{\bar{u}}),\bar{u}(\tau^{\bar{u}}))]\bar{h}(v(t),t)}{h^{\bar{u}}(\tau^{\bar{u}})}\bigg{]}\text{d}t
+\mathbb{E}\bigg{[} \Psi_x(\bar{X}(\tau^{\bar{u}}))^{\top}y(\tau^{\bar{u}})\big{)}\bigg{]}\\
&+\displaystyle\int_{0}^{\tau^{\bar{u}}}\mathbb{E}\bigg{[}
f_x(\bar{X}{(t)},\bar{u}(t))^{\top}y(t)+f_u(\bar{X}{(t)},\bar{u}(t))^{\top}v(t)  \bigg{]}\text{d}t+o(1).
\end{array}
\end{equation*}

{ Second, we consider the case $(ii)$. Notice that $\inf\bigg{\{}t:\mathbb{E}[\Phi(\bar{X}(t))]\geq \alpha,\ t\in [0,T] \bigg{\}}=T$, thus for sufficiently small $\rho$, if
$$
\tau^{{u_{\rho}}}=\inf\bigg{\{}t:\mathbb{E}[\Phi(X^{{u_{\rho}}}(t))]\geq \alpha,\ t\in [0,T] \bigg{\}}<T.
$$
Similar with the proof of case $(i)$, by case $(ii)$ of Lemma \ref{le-2}, one obtains
\begin{equation*}
\begin{array}
[c]{rl}
& \rho^{-1}\left[J(u_{\rho}(\cdot))-J(\bar{u}(\cdot))\right] \\
= &\displaystyle -\int_0^{\tau^{\bar{u}}}\bigg{[}\frac{ \tilde{\Psi}^{\bar{u}}(\tau^{\bar{u}})\bar{h}(v(t),t)}{h^{\bar{u}}(\tau^{\bar{u}})}+\frac{ \mathbb{E}[f(\bar{X}(\tau^{\bar{u}}),\bar{u}(\tau^{\bar{u}}))]\bar{h}(v(t),t)}{h^{\bar{u}}(\tau^{\bar{u}})}\bigg{]}\text{d}t
+\mathbb{E}\bigg{[} \Psi_x(\bar{X}(\tau^{\bar{u}}))^{\top}y(\tau^{\bar{u}})\big{)}\bigg{]}\\
&+\displaystyle\int_{0}^{\tau^{\bar{u}}}\mathbb{E}\bigg{[}
f_x(\bar{X}{(t)},\bar{u}(t))^{\top}y(t)+f_u(\bar{X}{(t)},\bar{u}(t))^{\top}v(t)  \bigg{]}\text{d}t+o(1).
\end{array}
\end{equation*}
If
$$
\inf\bigg{\{}t:\mathbb{E}[\Phi(X^{{u_{\rho}}}(t))]\geq \alpha,\ t\in [0,T] \bigg{\}}=+\infty,
$$
then, $\tau^{{u_{\rho}}}=T$, we can obtain
\begin{equation*}
\begin{array}
[c]{rl}
& \rho^{-1}\left[J(u_{\rho}(\cdot))-J(\bar{u}(\cdot))\right] \\
= &\displaystyle \mathbb{E}\bigg{[} \Psi_x(\bar{X}(\tau^{\bar{u}}))^{\top}y(\tau^{\bar{u}})\big{)}\bigg{]}+\displaystyle\int_{0}^{\tau^{\bar{u}}}\mathbb{E}\bigg{[}
f_x(\bar{X}{(t)},\bar{u}(t))^{\top}y(t)+f_u(\bar{X}{(t)},\bar{u}(t))^{\top}v(t)  \bigg{]}\mathrm{d}t+o(1).
\end{array}
\end{equation*}
Thus, the case $(ii)$ is right.

In the end, we consider case $(iii)$. Notice that $\bigg{\{}t:\mathbb{E}[\Phi(\bar{X}(t))]\geq \alpha,\ t\in [0,T] \bigg{\}}=\varnothing$, thus for sufficiently small $\rho$,
$$
\inf\bigg{\{}t:\mathbb{E}[\Phi(X^{{u_{\rho}}}(t))]\geq \alpha,\ t\in [0,T] \bigg{\}}=+\infty,
$$
and $\tau^{{u_{\rho}}}=T$, by case $(iii)$ of Lemma \ref{le-2}, we have
\begin{equation*}
\begin{array}
[c]{rl}
& \rho^{-1}\left[J(u_{\rho}(\cdot))-J(\bar{u}(\cdot))\right] \\
= &\displaystyle \mathbb{E}\bigg{[} \Psi_x(\bar{X}(\tau^{\bar{u}}))^{\top}y(\tau^{\bar{u}})\big{)}\bigg{]}+\displaystyle\int_{0}^{\tau^{\bar{u}}}\mathbb{E}\bigg{[}
f_x(\bar{X}{(t)},\bar{u}(t))^{\top}y(t)+f_u(\bar{X}{(t)},\bar{u}(t))^{\top}v(t)  \bigg{]}\mathrm{d}t+o(1).
\end{array}
\end{equation*}
} This completes the proof. $\ \ \ \ \ \ \ \ \Box$

\subsection{Stochastic maximum principle}\label{smp}
We introduce the following first-order  adjoint equation:
\begin{equation}%
\begin{array}
[c]{rl}%
-\text{d}{p}(t)= & \bigg{[}b_x(\bar{X}{(t)},\bar{u}(t))^{\top}p(t)+ \displaystyle\sum_{j=1}^d\sigma_x^j(\bar{X}{(t)},\bar{u}(t))^{\top}q^j(t) \\
               &-f_x(\bar{X}{(t)},\bar{u}(t))\bigg{]}\text{d}t-q(t)\text{d}W(t),\ t\in[0,\tau^{\bar{u}}),\\
p(\tau^{\bar{u}})= &-\Psi_{x}(\bar{X}(\tau^{\bar{u}}))^{\top}.
\end{array}
\label{prin-1}%
\end{equation}
{ Equation (\ref{prin-1}) is a linear Backward stochastic differential equation, we can obtain an explicit solution for equation (\ref{prin-1}) by the dual method, see Chapter 7 in \cite{Y99} for the basic theory of Backward stochastic differential equation}. Denoting
\begin{equation*}
{H}(x,u,p,q)=b(x,u)^{\top}p+\sum_{j=1}^d\sigma^j (x,u)^{\top}q^j-f(x,u),\text{ \  \ }%
(x,u,p,q)\in \mathbb{R}^m\times U\times \mathbb{R}^m\times \mathbb{R}^{m\times d}.%
\end{equation*}
The main result of this study is given as follows:
\begin{theorem}
\label{Max} Let Assumptions \ref{ass-b}, \ref{assb-b2} and \ref{ass-fai} hold,
$(\bar{u}(\cdot),\bar{X}(\cdot))$ be an optimal pair of (\ref{cost-2}). Suppose $h^{\bar{u}}(\tau^{\bar{u}})\neq 0$ and $h^{\bar{u}}(\cdot)$ is continuous at the point $\tau^{\bar{u}}$. Then, there exists  $(p(\cdot),q(\cdot))$ satisfying the series of first-order adjoint equations (\ref{prin-1}) and  the following hold.

\noindent (i). If $\tau^{\bar{u}}<T$, one obtains
\begin{equation}%
\begin{array}
[c]{ll}%
&H_u(\bar{X}(t),\bar{u}(t),p(t),q(t))(u-\bar{u}(t))+\displaystyle\frac{ \tilde{\Psi}^{\bar{u}}(\tau^{\bar{u}})\bar{h}(u-\bar{u}(t),t)}{h^{\bar{u}}(\tau^{\bar{u}})}+\frac{ \mathbb{E}[f(\bar{X}(\tau^{\bar{u}}),\bar{u}(\tau^{\bar{u}}))]\bar{h}(u-\bar{u}(t),t)}{h^{\bar{u}}
(\tau^{\bar{u}})}\leq 0,\\
\end{array}
\label{prin-2}%
\end{equation}
where
$$
\tilde{\Psi}^{\bar{u}}(\tau^{\bar{u}})=\mathbb{E}\bigg{[}\Psi_x(\bar{X}(\tau^{\bar{u}}))^{\top}b(\bar{X}(\tau^{\bar{u}}),\bar{u}(\tau^{\bar{u}}))
+\frac{1}{2}\sum_{j=1}^d\sigma^j(\bar{X}(\tau^{\bar{u}}),\bar{u}(\tau^{\bar{u}}))^{\top}\Psi_{xx}(\bar{X}(\tau^{\bar{u}}))
\sigma^j(\bar{X}(\tau^{\bar{u}}),\bar{u}(\tau^{\bar{u}}))\bigg{]},
$$
for any $u\in U$, $\text{a.e.}\ t \in[0,\tau^{\bar{u}})$, $P-\text{a.s.}$

\noindent(ii). If $\inf\bigg{\{}t:\mathbb{E}[\Phi(\bar{X}(t))]\geq \alpha,\ t\in [0,T] \bigg{\}}=T$, one obtains
\begin{equation}%
\begin{array}
[c]{ll}%
&H_u(\bar{X}(t),\bar{u}(t),p(t),q(t))(u-\bar{u}(t))+\displaystyle\frac{ \tilde{\Psi}^{\bar{u}}(\tau^{\bar{u}})\bar{h}(u-\bar{u}(t),t)}{h^{\bar{u}}(\tau^{\bar{u}})}+\frac{ \mathbb{E}[f(\bar{X}(\tau^{\bar{u}}),\bar{u}(\tau^{\bar{u}}))]\bar{h}(u-\bar{u}(t),t)}{h^{\bar{u}}
(\tau^{\bar{u}})}\leq 0,\\
\end{array}
\end{equation}
or
\begin{equation}%
\begin{array}
[c]{ll}%
&H_u(\bar{X}(t),\bar{u}(t),p(t),q(t))(u-\bar{u}(t))\leq 0,\\
\end{array}
\end{equation}
for any $u\in U$, $\text{a.e.}\ t \in[0,\tau^{\bar{u}})$, $P-\text{a.s.}$

\noindent(iii). If $\bigg{\{}t:\mathbb{E}[\Phi(\bar{X}(t))]\geq \alpha,\ t\in [0,T] \bigg{\}}=\varnothing$, we have
\begin{equation}%
\begin{array}
[c]{ll}%
&H_u(\bar{X}(t),\bar{u}(t),p(t),q(t))(u-\bar{u}(t))\leq 0,\\
\end{array}
\end{equation}
for any $u\in U$, $\text{a.e.}\ t \in[0,\tau^{\bar{u}})$, $P-\text{a.s.}$
\end{theorem}

\bigskip
\begin{remark}
In Theorem \ref{Max}, we propose a stochastic maximum principle for the varying terminal time optimal control problem. In case $(i)$, $\tau^{\bar{u}}<T$, there are two new terms:
$$
\displaystyle\frac{ \tilde{\Psi}^{\bar{u}}(\tau^{\bar{u}})\bar{h}(u-\bar{u}(t),t)}{h^{\bar{u}}(\tau^{\bar{u}})},\ \
\displaystyle
\frac{ \mathbb{E}[f(\bar{X}(\tau^{\bar{u}}),\bar{u}(\tau^{\bar{u}}))]\bar{h}(u-\bar{u}(t),t)}{h^{\bar{u}}
(\tau^{\bar{u}})},\ t\in[0,\tau^{\bar{u}}],
$$
where the first term is derived by $\Psi$, while the second term is derived by $f$. In fact, these two terms can be viewed as the penalty terms for the varying terminal time $\tau^{\bar{u}}$.
\end{remark}

In the following, we compare our new optimal control problem with the traditional one with and without state constraints. The traditional optimal control problem under state constraints is given as follows:

\begin{equation}
J(u(\cdot))=%
\mathbb{E}\bigg{[}{\displaystyle \int \limits_{0}^{T}}
f(X{(t)},u(t))\text{d}t+\Psi(X(T))\bigg{]},\label{ccost-2}%
\end{equation}
under state constraints
\begin{equation}
\label{cons-1}
\mathbb{E}\big{[}\Phi(X^u(T))\big{]}\geq \alpha.
\end{equation}

{ The following results can be found in Theorem 6.1 Chapter 3 of \cite{Y99}}.
\begin{theorem}
\label{c-th}
Let Assumptions \ref{ass-b}, \ref{assb-b2} and \ref{ass-fai} hold,
and let $(\hat{u}(\cdot),\hat{X}(\cdot))$ be an optimal pair of (\ref{ccost-2}) under state constraints (\ref{cons-1}).
Then, there exists $(\beta_0,\beta_1)\in \mathbb{R}^{2}$ satisfying
$$
\beta_0\geq 0,\ \ \left| \beta_0\right|^2+\left| \beta_1\right|^2=1,
$$
and
$$
\beta_1\big{(}\gamma-\mathbb{E}[\Phi(\hat{X}(T))]\big{)}\geq 0,\  \gamma\geq \alpha.
$$
The adapted solution $(p(\cdot),q(\cdot))$ satisfies the following  first-order adjoint equation:
\begin{equation}%
\begin{array}
[r]{rl}%
-\mathrm{d}{p}(t)=& \displaystyle \bigg{[}b_x(\hat{X}{(t)},\hat{u}(t))^{\top}p(t)+ \sum_{j=1}^d\sigma_x^j(\hat{X}{(t)},\hat{u}(t))^{\top}q^j(t) \\
               &-\beta_0f_x(\hat{X}{(t)},\hat{u}(t))\bigg{]}\mathrm{d}t-q(t)\mathrm{d}W(t),\ t\in[0,T),\\
p(T)=& -\beta_1\Psi_{x}(\hat{X}(T))^{\top},
\end{array}
\label{cprin-1}%
\end{equation}
and
\begin{equation}%
\begin{array}
[c]{ll}%
&H_u(\beta^0,\hat{X}(t),\hat{u}(t),p(t),q(t))(u-\hat{u}(t))
\leq 0,
\end{array}
\label{cprinc-2}%
\end{equation}
for any $u\in U$ and $t \in[0,T)$, where
\begin{equation*}
{H}(\beta_0,x,u,p,q)=b(x,u)^{\top}p+\sum_{j=1}^d\sigma^j (x,u)^{\top}q^j-\beta_0f(x,u),\text{ \  \ }%
(x,u,p,q)\in \mathbb{R}^m\times U\times \mathbb{R}^m\times \mathbb{R}^{m\times d}.%
\end{equation*}
\end{theorem}

We have established three kinds of maximum principle in Theorem \ref{Max}: In cases $(i)$ and $(ii)$, from Remark \ref{re-1}, we have $\mathbb{E}[\Phi(\bar{X}(\tau^{\bar{u}}))]=\alpha$ and $ \mathbb{E}[\Phi(\bar{X}(t))]<\alpha,\ t\in [0,\tau^{\bar{u}})$ for  an optimal pair $(\bar{u}(\cdot),\bar{X}(\cdot))$, which is more practical than the state constraints for the optimal control problem with $\mathbb{E}[\Phi(\bar{X}(T))]\geq \alpha$. However, we have $\mathbb{E}[\Phi(\bar{X}(\tau^{\bar{u}}))]=\alpha$ for an optimal pair $(\bar{u}(\cdot),\bar{X}(\cdot))$ of our new optimal control problem, while we have $\mathbb{E}[\Phi(\hat{X}(T))]\geq \alpha$ for an optimal pair $(\hat{u}(\cdot),\hat{X}(\cdot))$ of the traditional optimal control problem under state constraints:
 \begin{itemize}
\item { Notice that in Theorem \ref{c-th}, the parameter $(\beta_0,\beta_1)$ depends on the optimal pair $(\hat{u}(\cdot),\hat{X}(\cdot))$. We need to calculate the parameters  $(\beta_0,\beta_1)$ by the optimal pair $(\hat{u}(\cdot),\hat{X}(\cdot))$.  Thus, when we want to use Theorem \ref{c-th} to find an optimal pair  $(\hat{u}(\cdot),\hat{X}(\cdot))$ for cost functional (\ref{ccost-2}) under the state constraints (\ref{cons-1}), it is not easily to determine the parameter $(\beta_0,\beta_1)$.}
\item Our new optimal control model is based on the varying terminal time (\ref{time-1}), and involves minimizing the cost functional (\ref{cost-1}) with the varying terminal time. The constrained condition is introduced in the definition of the varying terminal time. The advantage of our model is that we can calculate the variation of the varying terminal time and obtain a stochastic maximum principle which can be easily verified.
\end{itemize}
In case $(iii), \ \bigg{\{}t:\mathbb{E}[\Phi(\bar{X}(t))]\geq \alpha,\ t\in [0,T] \bigg{\}}=\varnothing$,  which shows that $\mathbb{E}[\Phi(\bar{X}(t))]< \alpha,\ t\in [0,T]$ . Thus, our optimal control problem is same with the traditional optimal control problem without state constraints.

{In the following, we present an example to illustrate the application of Theorem \ref{Max} and compare with the traditional optimal control problem under state constraints}.
\begin{example}
Let $m=d=1,T=1,\ U=[1,2],\ \alpha=1$, $b(x,u)=x+u,\ \sigma(x,u)=0,\ f(x,u)=u, \ \Phi(x)=x$ and $\Psi(x)=0$. The controlled ordinal differential equation is given as follows:
\begin{equation}
X^u(s)=\int_0^s\bigg{[}X^u(t)+u(t)\bigg{]}\mathrm{d}t,
\end{equation}
where the varying terminal time is
$$
 \tau^{{u}}=\inf\bigg{\{}t:{X}^u(t)\geq \alpha,\ t\in [0,T] \bigg{\}}\bigwedge T,
$$
and the cost functional is
$$
J(u(\cdot))=\int_0^{\tau^{{u}}}{u}(s)\mathrm{d}s.
$$

We employ Theorem \ref{Max} to find an optimal pair $(\bar{u}(t),\bar{X}(t))$, $t\in [0,T]$. First, we suppose that $\tau^{\bar{u}}<1$. For $t\in [0,\tau^{\bar{u}}]$, we have $h^{\bar{u}}(t)=\bar{X}(t)+\bar{u}(t)$ and
$\bar{h}(v(t),t)=y(t)+v(t)$, where
\begin{equation*}%
\begin{array}
[c]{rl}%
\mathrm{d}{y}(t)= & \big{[}y(t)+v(t)\big{]}\mathrm{d}t, \ t\in[0,\tau^{\bar{u}}),\\
y(0)= &0.
\end{array}
\end{equation*}
Thus
\begin{equation}
\label{ex1}
\bar{h}(v(t),t)= \int_0^te^{t-s}v(s)\mathrm{d}s+v(t).
\end{equation}
The first-order  adjoint equation is
\begin{equation*}%
\begin{array}
[c]{rl}%
-\mathrm{d}{p}(t)= & p(t)\mathrm{d}t, \ t\in[0,\tau^{\bar{u}}),\\
p(\tau^{\bar{u}})= &0,
\end{array}
\end{equation*}
which implies that
\begin{equation}
\label{ex2}
p(t)=0,\ t\in [0,\tau^{\bar{u}}].
\end{equation}
By  $(i)$ of Theorem \ref{Max}, combining equations (\ref{ex1}), (\ref{ex2}) and $\bar{X}(\tau^{\bar{u}})=\alpha=1$  we have
\begin{equation*}%
\begin{array}
[c]{rl}%
&H_u(\bar{X}(t),\bar{u}(t),p(t),q(t))(u-\bar{u}(t))+\displaystyle\frac{ \tilde{\Psi}^{\bar{u}}(\tau^{\bar{u}})\bar{h}(u-\bar{u}(t),t)}{h^{\bar{u}}(\tau^{\bar{u}})}+
\frac{ f(\bar{X}(\tau^{\bar{u}}),\bar{u}(\tau^{\bar{u}}))\bar{h}(u-\bar{u}(t),t)}{h^{\bar{u}}
(\tau^{\bar{u}})}\\
=&(p(t)-1)(u-\bar{u}(t))+\displaystyle 0+\frac{\bar{u}(\tau^{\bar{u}})}
{1+\bar{u}(\tau^{\bar{u}})}\bar{h}(u-\bar{u}(t),t)\\
=&-(u-\bar{u}(t))\displaystyle+\frac{\bar{u}(\tau^{\bar{u}})}
{1+\bar{u}(\tau^{\bar{u}})}\bigg{[}\int_0^te^{t-s}(u-\bar{u}(s))\mathrm{d}s+u-\bar{u}(t) \bigg{]}\\
=&\displaystyle \bigg{[} \frac{\bar{u}(\tau^{\bar{u}})}
{1+\bar{u}(\tau^{\bar{u}})}-1 \bigg{]} (u-\bar{u}(t))+ \frac{\bar{u}(\tau^{\bar{u}})}
{1+\bar{u}(\tau^{\bar{u}})}\int_0^te^{t-s}(u-\bar{u}(s))\mathrm{d}s\\
\leq & 0.
\end{array}
\end{equation*}
It follows that
$$
\bigg{[}\frac{\bar{u}(\tau^{\bar{u}})}
{1+\bar{u}(\tau^{\bar{u}})} e^t-1 \bigg{]}u\leq \displaystyle \bigg{[} \frac{\bar{u}(\tau^{\bar{u}})}
{1+\bar{u}(\tau^{\bar{u}})}-1 \bigg{]} \bar{u}(t)+ \frac{\bar{u}(\tau^{\bar{u}})}
{1+\bar{u}(\tau^{\bar{u}})}\int_0^te^{t-s}\bar{u}(s)\mathrm{d}s,\ u\in [1,2],\ t\in [0,\tau^{\bar{u}}).
$$
Now, suppose that $\displaystyle \frac{\bar{u}(\tau^{\bar{u}})}
{1+\bar{u}(\tau^{\bar{u}})} e^{
\tau^{\bar{u}}}-1 \leq 0$. Then, we can obtain an optimal pair
\begin{equation}
\displaystyle   \bar{u}(t)=1,\ t\in [0,\tau^{\bar{u}}],
  \end{equation}
and
$$
\bar{X}(t)=e^t-1,\ t\in [0,\tau^{\bar{u}}].
$$
From $\bar{X}(\tau^{\bar{u}})=e^{\tau^{\bar{u}}}-1=1$, we have $\tau^{\bar{u}}=\ln 2<T=1$. Notice that
 $$
 \frac{\bar{u}(\tau^{\bar{u}})}
{1+\bar{u}(\tau^{\bar{u}})} e^{
\tau^{\bar{u}}}=0.5e^{
\tau^{\bar{u}}}\leq 1.
 $$
Thus, $\tau^{\bar{u}}$ satisfies the inequality $\displaystyle \frac{\bar{u}(\tau^{\bar{u}})}
{1+\bar{u}(\tau^{\bar{u}})} e^{
\tau^{\bar{u}}}-1 \leq 0$, which implies
$$
(\bar{u}(t),\bar{X}(t))=(1,e^t-1),\ t\in [ 0, \tau^{\bar{u}}],
$$
where $\tau^{\bar{u}}=\ln 2$.

{ Now, we consider the traditional optimal control problem under state constraints. The cost functional is given as follows:
$$
J_1(u(\cdot))=\int_0^{1}{u}(s)\mathrm{d}s,
$$
under constrained condition
$$
X^u(T)\geq 1.
$$
We can obtain that the optimal control is $\bar{u}_1(t)=1$ and the related optimal state is $\bar{X}_1(t)=e^t-1,\ t\in [0,1]$.

Notice that,
$$
(\bar{u}_1(t),\bar{X}_1(t))=(\bar{u}(t),\bar{X}(t))=(1,e^t-1),\ t\in [0, \tau^{\bar{u}}].
$$
In our varying terminal time optimal control problem, the cost functional is
$$
J(u(\cdot))=\int_0^{\tau^{u}}{u}(s)\mathrm{d}s,
$$
and the optimal terminal time is $\tau^{\bar{u}}=\ln 2<1$. These results indicate that we can stop to control the system at the optimal terminal time $\tau^{\bar{u}}$ which satisfies the constrained condition
$\bar{X}(\tau^{\bar{u}})=\bar{X}(\ln2)\geq 1$ and obtains a smaller cost functional, $J(\bar{u}(\cdot))=\ln2<1=J_1(\bar{u}_1(\cdot))$.

}
\end{example}

\bigskip

\noindent \textbf{The proof of Theorem \ref{Max}.} We first prove case $(i)$. For$\ t\in (
0,\tau^{\bar{u}}),$ applying It\^{o} formula  to $p(t)^{\top}y(t)$ gives%
\begin{equation*}%
\begin{array}
[c]{rl}
\mathrm{d}\left[p(t)^{\top}y(t)\right]=&\mathrm{d}\left[p(t)^{\top}\right]y(t)+p(t)^{\top}\mathrm{d}\left[y(t)\right]
+\mathrm{d}\left[p(t)^{\top}\right]\mathrm{d}\left[y(t)\right]\\
= &-\bigg{[}p(t)^{\top}b_x(\bar{X}{(t)},\bar{u}(t))+ \displaystyle\sum_{j=1}^d q^j(t)^{\top} \sigma_x^j(\bar{X}{(t)},\bar{u}(t))\\
               &-f_x(\bar{X}{(t)},\bar{u}(t))^{\top}\bigg{]}y(t)\mathrm{d}t+q(t)^{\top}y(t)\mathrm{d}W(t)\\
               &+p(t)^{\top}\bigg{[}b_{x}(\bar{X}{(t)},\bar{u}(t))y(t)+b_{u}(\bar{X}{(t)},\bar{u}(t)) v(t)\bigg{]}\mathrm{d}t\\
               &+p(t)^{\top} \displaystyle \sum_{j=1}^d\bigg{[} \sigma_{x}^j(\bar{X}{(t)},\bar{u}(t))y(t)+\sigma^j_u(\bar{X}{(t)},\bar{u}(t)) v(t)\bigg{]}\mathrm{d}W^j(t)\\
               &+\displaystyle \sum_{j=1}^d\bigg{[} q^j(t)^{\top}\sigma_{x}^j(\bar{X}{(t)},\bar{u}(t))y(t)+q^j(t)^{\top}\sigma^j_u(\bar{X}{(t)},\bar{u}(t)) v(t)\bigg{]}\mathrm{d}t.
\end{array}
\end{equation*}
Integrating on both sides of the above equation from $0$ to $\tau^{\bar{u}}$ and taking the expectation, one obtains
\begin{equation}%
\begin{array}
[c]{rl}
&\mathbb{E}\bigg{[} p(\tau^{\bar{u}})^{\top}y(\tau^{\bar{u}})-p(0)^{\top}y(0)\bigg{]}\\
=&\mathbb{E} \bigg{[}-\Psi_{x}(\bar{X}(\tau^{\bar{u}}))^{\top}y(\tau^{\bar{u}})\bigg{]}\\
=& \mathbb{E}\displaystyle \int \limits_{0}^{\tau^{\bar{u}}}\bigg{[}p(t)^{\top}b_u(\bar{X}(t),\bar{u}(t))v(t)+\sum_{j=1}^d
q^j(t)^{\top}\sigma_u^j (\bar{X}(t),\bar{u}(t))v(t)+f_x(\bar{X}(t),\bar{u}(t))^{\top}y(t)   \bigg{]}\mathrm{d}t.
\end{array}
\label{3max-1}%
\end{equation}
It follows that
\begin{equation*}%
\begin{array}
[c]{rl}
&\displaystyle  \mathbb{E} \bigg{[}-\Psi_{x}(\bar{X}(\tau^{\bar{u}}))^{\top}y(\tau^{\bar{u}})-\displaystyle \int \limits_{0}^{\tau^{\bar{u}}}
\big{(}f_{x}(\bar{X}{(t)},\bar{u}(t))^{\top}y(t)
+f_u(\bar{X}{(t)},\bar{u}(t))^{\top}v(t)\big{)}\text{d}t\bigg{]}\\
=&\displaystyle \displaystyle \int \limits_{0}^{\tau^{\bar{u}}}\mathbb{E}\bigg{[}p(t)^{\top}b_u(\bar{X}(t),\bar{u}(t))v(t)+\sum_{j=1}^d
q^j(t)^{\top}\sigma_u^j (\bar{X}(t),\bar{u}(t))v(t)
-f_u(\bar{X}(t),\bar{u}(t))^{\top}v(t)\bigg{]}\text{d}t.\\
\end{array}
\end{equation*}
{ By the definition of the function $H(\cdot)$ and $(i)$ of Lemma \ref{le-3}}, we have
\begin{equation*}%
\begin{array}
[c]{rl}
&\displaystyle\mathbb{E} \int \limits_{0}^{\tau^{\bar{u}}}
 H_u(\bar{X}(t),\bar{u}(t),p(t),q(t))v(t)\text{d}t\\
=&\displaystyle \mathbb{E} \int \limits_{0}^{\tau^{\bar{u}}}\bigg{[}p(t)^{\top}b_u(\bar{X}(t),\bar{u}(t))v(t)
+\sum_{j=1}^d
q^j(t)^{\top}\sigma_u^j (\bar{X}(t),\bar{u}(t))v(t)
-f_u(\bar{X}(t),\bar{u}(t))^{\top}v(t)\bigg{]}\mathrm{d}t\\
=& -\displaystyle \mathbb{E} \int_0^{\tau^{\bar{u}}}\bigg{[}\frac{ \tilde{\Psi}^{\rho}(\tau^{\bar{u}})\bar{h}(v(t),t)}{h^{\bar{u}}(\tau^{\bar{u}})}+\frac{ \mathbb{E}[f(\bar{X}(\tau^{\bar{u}}),\bar{u}(\tau^{\bar{u}}))]\bar{h}(v(t),t)}{h^{\bar{u}}(\tau^{\bar{u}})}\bigg{]}
\text{d}t
-\rho^{-1}\left[J(u_{\rho}(\cdot))-J(\bar{u}(\cdot))\right]+o(1)\\
\leq &-\displaystyle \mathbb{E}\int_0^{\tau^{\bar{u}}}\bigg{[}\frac{ \tilde{\Psi}^{\rho}(\tau^{\bar{u}})\bar{h}(v(t),t)}{h^{\bar{u}}(\tau^{\bar{u}})}+\frac{ \mathbb{E}[f(\bar{X}(\tau^{\bar{u}}),\bar{u}(\tau^{\bar{u}}))]\bar{h}(v(t),t)}{h^{\bar{u}}(\tau^{\bar{u}})}\bigg{]}
\text{d}t+o(1),
\end{array}
\end{equation*}
which implies that
$$
\displaystyle\mathbb{E} \int_0^{\tau^{\bar{u}}}\bigg{[} H_u(\bar{X}(t),\bar{u}(t),p(t),q(t))v(t)+\frac{ \tilde{\Psi}^{\rho}(\tau^{\bar{u}})\bar{h}(v(t),t)}{h^{\bar{u}}(\tau^{\bar{u}})}+\frac{ \mathbb{E}[f(\bar{X}(\tau^{\bar{u}}),\bar{u}(\tau^{\bar{u}}))]\bar{h}(v(t),t)}{h^{\bar{u}}
(\tau^{\bar{u}})}\bigg{]}\text{d}t\leq o(1).
$$

Notice that  for any $u\in U$, we set $v(\cdot)=u-\bar{u}(\cdot)$. Thus, $\bar{u}(\cdot)+v(\cdot)\in \mathcal{U}[0,T]$. Letting $\rho\to 0$, we obtain
\begin{equation*}%
\begin{array}
[c]{ll}%
&H_u(\bar{X}(t),\bar{u}(t),p(t),q(t))(u-\bar{u}(t))+\displaystyle\frac{ \tilde{\Psi}^{\rho}(\tau^{\bar{u}})\bar{h}(u-\bar{u}(t),t)}{h^{\bar{u}}(\tau^{\bar{u}})}+\frac{ \mathbb{E}[f(\bar{X}(\tau^{\bar{u}}),\bar{u}(\tau^{\bar{u}}))]\bar{h}(u-\bar{u}(t),t)}{h^{\bar{u}}
(\tau^{\bar{u}})}\leq 0,\\
\end{array}
\end{equation*}
for any $u\in U$,  $\text{a.e.}\ t \in[0,\tau^{\bar{u}})$, and $P-\text{a.s.}$ If not, we can prove the above results by contradiction.

{ Second, we consider the case $(ii)$. Notice that $\inf\bigg{\{}t:\mathbb{E}[\Phi(\bar{X}(t))]\geq \alpha,\ t\in [0,T] \bigg{\}}=T$, thus for sufficiently small $\rho$, if
$$
\tau^{{u_{\rho}}}=\inf\bigg{\{}t:\mathbb{E}[\Phi(X^{{u_{\rho}}}(t))]\geq \alpha,\ t\in [0,T] \bigg{\}}<T.
$$
By the case $(ii)$ of Lemma \ref{le-3}, we have
\begin{equation*}
\begin{array}
[c]{rl}
& \rho^{-1}\left[J(u_{\rho}(\cdot))-J(\bar{u}(\cdot))\right] \\
= &\displaystyle -\int_0^{\tau^{\bar{u}}}\bigg{[}\frac{ \tilde{\Psi}^{\bar{u}}(\tau^{\bar{u}})\bar{h}(v(t),t)}{h^{\bar{u}}(\tau^{\bar{u}})}+\frac{ \mathbb{E}[f(\bar{X}(\tau^{\bar{u}}),\bar{u}(\tau^{\bar{u}}))]\bar{h}(v(t),t)}{h^{\bar{u}}(\tau^{\bar{u}})}\bigg{]}\text{d}t
+\mathbb{E}\bigg{[} \Psi_x(\bar{X}(\tau^{\bar{u}}))^{\top}y(\tau^{\bar{u}})\big{)}\bigg{]}\\
&+\displaystyle\int_{0}^{\tau^{\bar{u}}}\mathbb{E}\bigg{[}
f_x(\bar{X}{(t)},\bar{u}(t))^{\top}y(t)+f_u(\bar{X}{(t)},\bar{u}(t))^{\top}v(t)  \bigg{]}\text{d}t+o(1).
\end{array}
\end{equation*}
Then, similar with the proof of case $(i)$, we can obtain
\begin{equation*}%
\begin{array}
[c]{ll}%
&H_u(\bar{X}(t),\bar{u}(t),p(t),q(t))(u-\bar{u}(t))+\displaystyle\frac{ \tilde{\Psi}^{\rho}(\tau^{\bar{u}})\bar{h}(u-\bar{u}(t),t)}{h^{\bar{u}}(\tau^{\bar{u}})}+\frac{ \mathbb{E}[f(\bar{X}(\tau^{\bar{u}}),\bar{u}(\tau^{\bar{u}}))]\bar{h}(u-\bar{u}(t),t)}{h^{\bar{u}}
(\tau^{\bar{u}})}\leq 0.\\
\end{array}
\end{equation*}
If
$$
\inf\bigg{\{}t:\mathbb{E}[\Phi(X^{{u_{\rho}}}(t))]\geq \alpha,\ t\in [0,T] \bigg{\}}=+\infty,
$$
by the case $(ii)$ of Lemma \ref{le-3}, we have
\begin{equation*}
\begin{array}
[c]{rl}
& \rho^{-1}\left[J(u_{\rho}(\cdot))-J(\bar{u}(\cdot))\right] \\
= &\displaystyle \mathbb{E}\bigg{[} \Psi_x(\bar{X}(\tau^{\bar{u}}))^{\top}y(\tau^{\bar{u}})\big{)}\bigg{]}+\displaystyle\int_{0}^{\tau^{\bar{u}}}\mathbb{E}\bigg{[}
f_x(\bar{X}{(t)},\bar{u}(t))^{\top}y(t)+f_u(\bar{X}{(t)},\bar{u}(t))^{\top}v(t)  \bigg{]}\mathrm{d}t+o(1).
\end{array}
\end{equation*}
It follows that
\begin{equation*}%
\begin{array}
[c]{ll}%
&H_u(\bar{X}(t),\bar{u}(t),p(t),q(t))(u-\bar{u}(t))\leq 0.\\
\end{array}
\end{equation*}

In the end, we consider case $(iii)$. Notice that $\bigg{\{}t:\mathbb{E}[\Phi(\bar{X}(t))]\geq \alpha,\ t\in [0,T] \bigg{\}}=\varnothing$, thus for sufficiently small $\rho$,
$$
\inf\bigg{\{}t:\mathbb{E}[\Phi(X^{{u_{\rho}}}(t))]\geq \alpha,\ t\in [0,T] \bigg{\}}=+\infty,
$$
and $\tau^{{u_{\rho}}}=T$, by case $(iii)$ of Lemma \ref{le-3}, we have
\begin{equation*}
\begin{array}
[c]{rl}
& \rho^{-1}\left[J(u_{\rho}(\cdot))-J(\bar{u}(\cdot))\right]\\
= &\displaystyle \mathbb{E}\bigg{[} \Psi_x(\bar{X}(\tau^{\bar{u}}))^{\top}y(\tau^{\bar{u}})\big{)}\bigg{]}+\displaystyle\int_{0}^{\tau^{\bar{u}}}\mathbb{E}\bigg{[}
f_x(\bar{X}{(t)},\bar{u}(t))^{\top}y(t)+f_u(\bar{X}{(t)},\bar{u}(t))^{\top}v(t)  \bigg{]}\mathrm{d}t+o(1).
\end{array}
\end{equation*}
Thus, we can obtain
\begin{equation*}%
\begin{array}
[c]{ll}%
&H_u(\bar{X}(t),\bar{u}(t),p(t),q(t))(u-\bar{u}(t))\leq 0.
\end{array}
\end{equation*}
} This completes the proof. $\ \ \ \ \ \ \ \ \Box$

\section{Conclusion}
To solve the stochastic optimal control problem under state constraints, we introduce a varying terminal time optimal control problem structure, in which we can simultaneously balance the terminal time and minimize the cost functional. We investigate three different cases maximum principles for this new optimal control problem, including the traditional optimal control problem as one of them. {In addition, we can view this optimal control problem as a new structure for the optimal control problem under state constraints, in which we define a varying terminal time via the constrained condition. Furthermore, employing the convex variation method, we establish a novel stochastic maximum principle for this new optimal control problem}. This paper presents the first step in considering this new optimal control problem, based on which we can continue to work on  topics such as the dynamic programming principle, the relationship between the stochastic maximum principle and the dynamic programming principle, and the stochastic linear quadratic optimal control problem.

\appendix
\section{General $\Phi(\cdot)$}\label{app-a}
In this section, we calculate the function $\bar{h}(v(\cdot),\cdot)$ for general $\Phi(\cdot)$.
Note that, for $t\in [0,\tau^{\bar{u}}]$,
\begin{equation}%
\begin{array}
[c]{ll}%
&h^{{u}}(t)=\mathbb{E}\bigg{[}\Phi_x(X^{{u}}(t))^{\top}b(X^{{u}}(t),{u}(t))
+\displaystyle \frac{1}{2}\sum_{j=1}^d\sigma^j(X^{{u}}(t),{u}(t))^{\top}\Phi_{xx}(X^{{u}}(t))\sigma^j(X^{{u}}(t)
,{u}(t))\bigg{]},
\end{array}
\end{equation}
and
\begin{equation}
\mathbb{E}[\Phi(X^u(t))]=\Phi(x_0)+\int_0^th^u(s)\text{d}s.
\end{equation}
The function $\bar{h}(v(\cdot),\cdot)$ is defined as
\begin{equation}
\begin{array}
[c]{rl}%
&\bar{h}(v(t),t)=\displaystyle \lim_{\rho\to 0}\frac{h^{u_{\rho}}(t)-h^{\bar{u}}(t)}{\rho},\
t\in [0,\tau^{\bar{u}}].\\
\end{array}
\end{equation}

For $t\in [0,\tau^{\bar{u}}]$, we set
$$
g(X^u(t),u(t))=\Phi_x(X^{{u}}(t))^{\top}b(X^{{u}}(t),{u}(t))
+\displaystyle \frac{1}{2}\sum_{j=1}^d\sigma^j(X^{{u}}(t),{u}(t))^{\top}\Phi_{xx}(X^{{u}}(t))\sigma^j(X^{{u}}(t)
,{u}(t)).
$$
By adjoint equation (\ref{apro-1}), it follows that
\begin{equation}
\begin{array}
[c]{rl}%
&\displaystyle \lim_{\rho\to 0}\frac{h^{u_{\rho}}(t)-h^{\bar{u}}(t)}{\rho} =\mathbb{E}\bigg{[}g_{x}(\bar{X}{(t)},\bar{u}(t))^{\top}y(t)
+g_{u}(\bar{X}{(t)},\bar{u}(t)) v(t)\bigg{]},
\end{array}
\end{equation}
where
\begin{equation}
\begin{array}
[c]{rl}%
&g_{x}(\bar{X}{(t)},\bar{u}(t))\\
=& \Phi_{xx}(\bar{X}{(t)})^{\top}b(\bar{X}{(t)},\bar{u}(t))+
b_x(\bar{X}{(t)},\bar{u}(t))^{\top}\Phi_x(\bar{X}{(t)})\\
&+\displaystyle \sum_{j=1}^d\sigma_x^j(\bar{X}{(t)},\bar{u}(t))^{\top}\Phi_{xx}(\bar{X}{(t)})\sigma^j(\bar{X}{(t)}
,\bar{u}(t))\\
&+\displaystyle \frac{1}{2}\sum_{j=1}^d\Phi_{xxx}(\bar{X}{(t)})\sigma^j(\bar{X}{(t)}
,\bar{u}(t))\sigma^j(\bar{X}{(t)},\bar{u}(t)),
\end{array}
\end{equation}
and
\begin{equation}
\begin{array}
[c]{rl}%
&g_{u}(\bar{X}{(t)},\bar{u}(t))
=b_u(\bar{X}{(t)},\bar{u}(t))^{\top} \Phi_x(\bar{X}{(t)})+\displaystyle \sum_{j=1}^d\sigma_u^j(\bar{X}{(t)},\bar{u}(t))^{\top}\Phi_{xx}(\bar{X}{(t)})\sigma^j(\bar{X}{(t)}
,\bar{u}(t)).
\end{array}
\end{equation}

Similar to the proof of the traditional stochastic maximum principle  with a convex control domain, we can introduce a new adjoint equation to derive a dual representation for
$$
\mathbb{E}\displaystyle \int_0^{\tau^{\bar{u}}}\bigg{[}g_{x}(\bar{X}{(t)},\bar{u}(t))^{\top}y(t)
+g_{u}(\bar{X}{(t)},\bar{u}(t)) v(t)\bigg{]}\mathrm{d}t.
$$
For a convex optimal control set $U$, we introduce the following first-order  adjoint equations:
\begin{equation}%
\begin{array}
[c]{rl}%
-\text{d}{p}_0(t)= & \bigg{[}b_x(\bar{X}{(t)},\bar{u}(t))^{\top}p_0(t)+ \displaystyle\sum_{j=1}^d\sigma_x^j(\bar{X}{(t)},\bar{u}(t))^{\top}q_0^j(t) \\
               &-g_x(\bar{X}{(t)},\bar{u}(t))\bigg{]}\text{d}t-q_0(t)\text{d}W(t),\ \ t\in[0,\tau^{\bar{u}})\\
p_0(\tau^{\bar{u}})= &0.
\end{array}
\label{aaprin-1}%
\end{equation}
Denoting
\begin{equation*}
\mathcal{H}(x,u,p,q)=b(x,u)^{\top}p+\sum_{j=1}^d\sigma^j (x,u)^{\top}q^j-g(x,u),\text{ \  \ }%
(x,u,p,q)\in \mathbb{R}^m\times U\times \mathbb{R}^m\times \mathbb{R}^{m\times d}.%
\end{equation*}
\begin{theorem}
\label{Max2} Let Assumptions \ref{ass-b}, \ref{assb-b2} and \ref{ass-fai} hold,
and $(\bar{u}(\cdot),\bar{X}(\cdot))$ be an optimal pair of (\ref{cost-2}). We have
\begin{equation*}%
\begin{array}
[c]{rl}%
&\mathbb{E}\displaystyle \int_0^{\tau^{\bar{u}}}\bigg{[}g_{x}(\bar{X}{(t)},\bar{u}(t))^{\top}y(t)
+g_{u}(\bar{X}{(t)},\bar{u}(t)) v(t)\bigg{]}\mathrm{d}t
=\displaystyle \mathbb{E}\int_0^{\tau^{\bar{u}}}\bigg{[} -\mathcal{H}_u(\bar{X}(t),\bar{u}(t),p_0(t),q_0(t))v(t)\bigg{]}\mathrm{d}t.\\
\end{array}
\end{equation*}
\end{theorem}
\noindent\textbf{Proof:} This proof is similar to that of Theorem \ref{Max}. For$\ t\in (
0,\tau^{\bar{u}}),$ applying It\^{o} formula  to $p_0(t)^{\top}y(t)$,
\begin{equation*}%
\begin{array}
[c]{rl}
\mathrm{d}\left[p_0(t)^{\top}y(t)\right]=&\mathrm{d}\left[p_0(t)^{\top}\right]y(t)
+p_0(t)^{\top}\mathrm{d}\left[y(t)\right]
+\mathrm{d}\left[p_0(t)^{\top}\right]\mathrm{d}\left[y(t)\right]\\
= &-\bigg{[}p_0(t)^{\top}b_x(\bar{X}{(t)},\bar{u}(t))+ \displaystyle\sum_{j=1}^d q_0^j(t)^{\top} \sigma_x^j(\bar{X}{(t)},\bar{u}(t))\\
               &-g_x(\bar{X}{(t)},\bar{u}(t))^{\top}\bigg{]}y(t)\mathrm{d}t
               +q_0(t)^{\top}y(t)\mathrm{d}W(t)\\
               &+p_0(t)^{\top}\bigg{[}b_{x}(\bar{X}{(t)},\bar{u}(t))y(t)
               +b_{u}(\bar{X}{(t)},\bar{u}(t)) v(t)\bigg{]}\mathrm{d}t\\
               &+p_0(t)^{\top} \displaystyle \sum_{j=1}^d\bigg{[} \sigma_{x}^j(\bar{X}{(t)},\bar{u}(t))y(t)+\sigma^j_u(\bar{X}{(t)},\bar{u}(t)) v(t)\bigg{]}\mathrm{d}W^j(t)\\
               &+\displaystyle \sum_{j=1}^d\bigg{[} q_0^j(t)^{\top}\sigma_{x}^j(\bar{X}{(t)},\bar{u}(t))y(t)+q_0^j(t)^{\top}\sigma^j_u(\bar{X}{(t)},\bar{u}(t)) v(t)\bigg{]}\mathrm{d}t.
\end{array}
\end{equation*}
Integrating  both sides of the above equation from $0$ to $\tau^{\bar{u}}$ and taking the expectation, one obtains
\begin{equation*}%
\begin{array}
[c]{rl}
&\mathbb{E}\bigg{[} p_0(\tau^{\bar{u}})^{\top}y(\tau^{\bar{u}})-p_0(0)^{\top}y(0)\bigg{]}\\
=& \mathbb{E}\displaystyle \int \limits_{0}^{\tau^{\bar{u}}}\bigg{[}p_0(t)^{\top}b_u(\bar{X}(t),\bar{u}(t))v(t)+\sum_{j=1}^d
q_0^j(t)^{\top}\sigma_u^j (\bar{X}(t),\bar{u}(t))v(t)+g_x(\bar{X}(t),\bar{u}(t))^{\top}y(t)   \bigg{]}\mathrm{d}t.
\end{array}
\end{equation*}
It follows that
\begin{equation*}%
\begin{array}
[c]{rl}
&\displaystyle  -\mathbb{E}\displaystyle \int \limits_{0}^{\tau^{\bar{u}}}
 \bigg{[}g_{x}(\bar{X}{(t)},\bar{u}(t))^{\top}y(t)
+g_u(\bar{X}{(t)},\bar{u}(t))^{\top}v(t)\bigg{]}\text{d}t\\
=&\displaystyle \mathbb{E}\int \limits_{0}^{\tau^{\bar{u}}}\bigg{[}p_0(t)^{\top}b_u(\bar{X}(t),\bar{u}(t))v(t)+\sum_{j=1}^d
q_0^j(t)^{\top}\sigma_u^j (\bar{X}(t),\bar{u}(t))v(t)
-g_u(\bar{X}(t),\bar{u}(t))^{\top}v(t)\bigg{]}\text{d}t.\\
\end{array}
\end{equation*}
By the definition of the  function $\mathcal{H}(\cdot)$, we have
\begin{equation*}%
\begin{array}
[c]{rl}
&\displaystyle \mathbb{E} \int \limits_{0}^{\tau^{\bar{u}}}
 \mathcal{H}_u(\bar{X}(t),\bar{u}(t),p_0(t),q_0(t))v(t)\text{d}t\\
=&\displaystyle \mathbb{E} \int \limits_{0}^{\tau^{\bar{u}}}\bigg{[}p_0(t)^{\top}b_u(\bar{X}(t),\bar{u}(t))v(t)
+\sum_{j=1}^d
q_0^j(t)^{\top}\sigma_u^j (\bar{X}(t),\bar{u}(t))v(t)
-g_u(\bar{X}(t),\bar{u}(t))^{\top}v(t)\bigg{]}\mathrm{d}t,\\
\end{array}
\end{equation*}
which implies that
\begin{equation*}%
\begin{array}
[c]{rl}%
&\mathbb{E}\displaystyle \int_0^{\tau^{\bar{u}}}\bigg{[}g_{x}(\bar{X}{(t)},\bar{u}(t))^{\top}y(t)
+g_{u}(\bar{X}{(t)},\bar{u}(t)) v(t)\bigg{]}\mathrm{d}t
=\displaystyle \mathbb{E}\int_0^{\tau^{\bar{u}}}\bigg{[} -\mathcal{H}_u(\bar{X}(t),\bar{u}(t),p_0(t),q_0(t))v(t)\bigg{]}\mathrm{d}t.\\
\end{array}
\end{equation*}
This completes the proof.$\ \ \ \ \ \ \ \ \ \Box$

\end{document}